\documentclass[a4paper,12pt,twoside]{article}
\usepackage{amssymb}
\usepackage[final]{lpretex}
\usepackage{title}
\usepackage{a4}
\usepackage{amscd}
\usepackage{amsmath}
\input xy
\xyoption{all}

\EndOfDump

\def\DHpartformat#1{(#1) }
\def\DHrefpart#1{(\DHRefpart{#1})}
\LBsetstyleof{partno}{arabic}

\DocVersion[Changed to LaTeX]{1}
\DocVersion[Added deformation theory, rearranged; Sep 09 1999]{2}
\DocVersion[Two sections 1 combined; Sep 16 1999]{3}
\DocVersion[Minor tweaks; Sep 17 1999]{4}
\DocVersion[Added a new section and minor fixes; Sep 23 1999]{5}
\DocVersion[Many small changes (NSB), One new section (TE); Wed Oct  6  1999]{6}
\DocVersion[Many small changes (NSB), Stacks (TE); Wed Oct  11  1999]{7}

\let\define\def

\def\A {{\mathbb A}}  \def\C {{\mathbb C}}
  \def\F {{\mathbb F}}

\def\GG {{\mathbb G}}   

\def\N {{\mathbb N}}  \def\P {{\mathbb P}} 
\def\Q {{\mathbb Q}} \def\R {{\mathbb R}}
 
 \def\W {{\mathbb W}} 
\def\Z {{\mathbb Z}} 

\define \n {\mathbb N}
\define \z {\mathbb Z}
\define \q {\mathbb Q}
\define \PP {\mathbb P}

\def\sA {{\Cal A}}  \def\sC {{\Cal C}}
\def\sD {{\Cal D}} \def\sE {{\Cal E}} \def\sF {{\Cal F}}
  \def\sI {{\Cal I}}
  
 \def\sN {{\Cal N}} \def\sO {{\Cal O}}
 
\def\sS {{\Cal S}}  
  \def\sX {{\Cal X}}

\def\sZ {{\Cal Z}}
\define \cN {\Cal N}
\define \cf {\Cal F}
\define \cg {\Cal G}
\define \cE {\Cal E}
\define \ce {\Cal E}
\define \cc {\Cal C}
\define \cV {\Cal V}
\define \cA {\Cal A}
\define \cK {\Cal K}
\define \cO {\Cal O}
\define \cF {\Cal F}
\define \cn {\Cal N}
\define \cI {\Cal I}
\define \sP {\Cal P}
\def\a {\alpha} 
\def\g {\gamma} 

\def\s {\sigma}  
 \def\L {\Lambda} 
 
\define\x {\xi}
\define\y {\eta}
\define\G {\Gamma}
\define\r {\rho}
\define\w {\omega}

\def\tX {\widetilde X}

\def \tD {\widetilde D}

\def \trho {\tilde {\rho}}

\def \tp {\widetilde{\mathbb P}}
\define \tH {\widetilde H}
\define \tG {\widetilde{\Gamma}}
\define \tW {\widetilde W}
\define \tF {\widetilde F}
\define \tm {\tilde m}
\define \St {\widetilde S}
\define \Xt {\widetilde X}
\define \tS {\widetilde S}
\define \tpsi {\tilde \psi}
\define \tL {\widetilde L}
\define \tE {\widetilde E}
\define \tl {\tilde l}
\define \tA {\widetilde A}
\define \tom {\tilde\omega}
\define \tT {\widetilde T}
\define \tB {\widetilde B}
\define \tf {\tilde f}
\define \tsA {\widetilde{\sA}}

\define \tM {\widetilde M}
\define \tphi {\widetilde{\phi}}
\define \trho {\widetilde{\rho}}
\define \tR {\widetilde R}
\define \tQ {\widetilde Q}
\define \tp {\tilde p}
\define \tq {\tilde q}
\define \tc {\tilde c}
\define \tsF {\widetilde {\sF}}
\define \tsO {\widetilde {\sO}}
\define \tx {\tilde x}
\define \tg {\tilde g}
\define \tw {\tilde w}

\def\pd {\partial}

\def \Dx1 {\frac{\pd}{{\pd} x_1}}
\def \Dy1 {\frac{\pd}{{\pd} y_1}}
\def \Dz1 {\frac{\pd}{{\pd} z_1}}
\def \Dx2 {\frac{\pd}{{\pd} x_2}}
\def \Dy2 {\frac{\pd}{{\pd} y_2}}
\def \Dz2 {\frac{\pd}{{\pd} z_2}}

\def\q {\quad} 

\def\mapdiagr#1{\Big\searrow\rlap{$\raise 5pt\vbox{{\hbox{$\mkern -15mu\scriptstyle#1$}}}$}}   

\def\mapdiagl#1{\llap{$\raise 5pt\vbox{{\hbox{$\scriptstyle#1\mkern
-15mu$}}}$}\Big\swarrow}              

\def\Mapdiagr#1{\nearrow\rlap{$\lower 5pt\vbox{{\hbox{$\mkern
-15mu\scriptstyle#1$}}}$}} 

\def\Mapdiagl#1{\llap{$\lower 5pt\vbox{{\hbox{$\scriptstyle#1\mkern
-15mu$}}}$}\searrow} 

\def\Mapswr#1{\swarrow\rlap{$\lower 5pt\vbox{{\hbox{$\mkern
-15mu\scriptstyle#1$}}}$}}              

\def\Mapnwl#1{\nwarrow\rlap{$\lower 5pt\vbox{{\hbox{$\mkern
-15mu\scriptstyle#1$}}}$}}

\def \inj {\hookrightarrow}

\define \Rhook {\hookrightarrow}

\def \half {\raise1pt\hbox{$\scriptstyle
        \frac{1}{2}\displaystyle$}}

\def \x{{\sl X}\llap{$\mkern -2mu {\scriptstyle -}$}}
\def \Hom {\operatorname{Hom}}

\def \Res {\operatorname{Res}}

\def \Pic {\operatorname{Pic}}
\def \Sing {\operatorname{Sing}}
\define \Kod {\operatorname{Kod}}
\define \dimension {\operatorname{dim}}
\define \codim {\operatorname{codim}}
\define \contr {\operatorname{contr}}
\define \rk {\operatorname{rank}}
\define \im {\operatorname{im}}
\define \Mor {\operatorname{Mor}}
\define \Cl {\operatorname{Cl}}
\define \Hilb {\operatorname{Hilb}}
\define \degree {\operatorname{deg}}
\define \mult {\operatorname{mult}}
\define \Aut {\operatorname{Aut}}
\define \NS {\operatorname{NS}}
\define \Gal {\operatorname{Gal}}
\define \ch {\operatorname{char}}
\define \Jac {\operatorname{Jac}}
\define \Km {\operatorname{Km}}
\define \Sec {\operatorname{Sec}}
\define \Stab {\operatorname{Stab}}
\define \Br {\operatorname{Br}}
\define \inv {\operatorname{inv}}
\define \tr {\operatorname{tr}}
\define \Frob {\operatorname{Frob}}
\define \Symn {\operatorname{Sym}^n}
\define \Ev {\sE^\vee}
\define \ordp {\operatorname{ord}_p}
\define \Supp {\operatorname{Supp}}
\define \Ann {\operatorname{Ann}}
\define \disc {\operatorname{disc}}
\define \Lie {\operatorname{Lie}}
\define \embdim {\operatorname{embdim}}

\def \Fix{\operatorname{Fix}}

\def\Tr{\operatorname{Tr}}

\def\hD{\widehat D}
\def\hb{\hat\beta}

\def\Def{\operatorname{Def}}
\def\NE{\operatorname{NE}}

\def\W{\cansymb{W}}

\def\hod#1#2#3#4{\ensuremath{\if#30 H^{#2}({#1},{\cal O}_{#1}) \else 
 H^{#2}(#1,\Omega^{#3}\if\relax{#4}\relax_{#1}\else _{#1/#4}\fi)\fi}}

\begin{document}
\title[Weyl group covers and the integral cohomology of $G/P$]
{Weyl group covers for
Brieskorn's resolutions
in all characteristics and the integral cohomology of $G/P$}

\author{N. I. Shepherd-Barron}
\address
{Math. Dept.\\
King's College London\\
Strand\\
London WC2R 2LS UK}
\email{Nicholas.Shepherd-Barron@kcl.ac.uk}

\maketitle
\begin{section}{Introduction}

Artin proved [Ar74] that, over
a field $k$ of any characteristic,
for a given affine surface $X_s$
over $k$ with only rational singularities
there is a unique irreducible component $A$
(``the Artin component'')
of the deformation space of $X_s$
that contains all those deformations that
can be simultaneously resolved after some finite
covering of the base. 
This extended earlier constructions by Brieskorn
that concerned rational double points 
(these are
also called du Val singularities,
Kleinian singularities, rational double points,
simple singularities, ...)
in characteristic zero.
Burns and Rapoport conjectured ([BR75], Conjecture 7.4)
that this covering of $A$ is Galois and that its Galois group 
is the Weyl
group $W$ corresponding to the configuration of $(-2)$-curves
in the minimal resolution $X'_s$ of $X_s$.

The main results of this paper are these.

\begin{theorem} (= Theorem \ref{main})
For any field $k$ the conjecture of Burns and Rapoport is true 
over $k$ and over any coefficient ring $\L$ for $k$.
\noproof
\end{theorem}

When restricted to the context of certain RDPs defined over $\Z$
Theorem \ref{main} can be stated and proved over $\Z$
and then we learn something
about the invariant theory of Weyl groups over $\Z$
and about the integral cohomology rings $H^*(G/P,\Z)$,
where $G/P$ is a partial flag variety, as follows.
 
Recall that
if $L$ is a weight lattice then
$W$ acts on the polynomial $\Z$-algebra $\Z[L]$ 
and, over a ring $R$ in which a certain integer $m$
is invertible (for example, $m=30$ in the case of $E_8$),
the corresponding ring of invariants $R[L]^{W}$
is a polynomial $R$-algebra [Dm73]. This is false
for $E_8$
over $\Z$ or over a field of characteristic $2,3$ or $5$.
However, we prove that, for the corresponding root lattice $M$
rather than for the weight lattice, it is stably true in the $ADE$ case, in that
it becomes true after ``polynomial enlargement''.

\begin{theorem} (= Theorem \ref{extend})
There is an effective graded action of $W$ 
on a graded polynomial $\Z[M]$-algebra
$\sO$ 
such that the ring $\sO^{W}$ of invariants
is a polynomial $\Z$-algebra.
\noproof
\end{theorem}
However, we do not know how to write down the $W$-action.
That is, we do not have formulae for the action of any
of the reflexions in $W$, not even the simple ones.
Nor do we have formulae that describe generators of $\sO^W$
in terms of the given generators of $\sO$.

Let $G$ denote a simple algebraic group of type $ADE$. 
Using the results and techniques due to Bernstein, Gel'fand and Gel'fand [BGG73]
and Demazure [Dm73] we are lead to a
description of the integral cohomology ring
$H^*(G/B,\Z)$
(which ring equals the Chow ring, up to a doubling of degrees,
since the flag variety $G/B$ is paved
by affine spaces) that is an integral version of Borel's
description of $H^*(G/B,\Q)$ as the ring of co-invariants
associated to the action of the Weyl group on the $\Q$-polynomial
algebra $\Q[M]$.

\begin{corollary} (= Corollary \ref{5.11})
$H^*(G/B,\Z)$ is isomorphic, as a graded ring
with an action of $W$, to the ring
$\sO/(\sO^W_+.\sO)$ of co-invariants.
\noproof
\end{corollary}

Via the results of [BGG73] this leads to a similar description of $H^*(G/P,\Z)$
where $P$ is a  parabolic subgroup $G$ that corresponds
to a given set $\Theta$ of simple roots.
Let $W_\Theta$ denote the corresponding subgroup of $W$.

\begin{corollary} (= Corollary \ref{BGG})
$\sO^{W_\Theta}$ is also a polynomial ring over $\Z$ and
$H^*(G/P,\Z)$ is isomorphic
to the quotient ring
$\sO^{W_\Theta}/(\sO^W_+.\sO^{W_\Theta})$.
\noproof
\end{corollary}
These results complement earlier
explicit descriptions of $H^*(G/B,\Z)$ in terms of generators
and relations. See, for example,  [TW74] for
$G=D_n$ and
[Na10] and [DZ15] for $G=E_8$.

We also explain, in Theorem \ref{enriques},
some consequences for the local moduli of Enriques surfaces
in characteristic two.

We now give 
a more detailed sketch of the background.
The conjecture of [BR75] was proved by Wahl [Wa79] in characteristic zero;
his proof depends upon the fact that this conjecture
had already been proved, if the characteristic is either zero
or bigger than the relevant Coxeter number, when $X_s$ is an RDP
in earlier work by
many people (Brieskorn [Br70], Tyurina [Ty70], Slodowy [Sl80]).

They showed that one way of achieving simultaneous resolution
for du Val singularities is to embed the picture into
the corresponding simply connected simple algebraic group.
This approach also gives a complete description
of the finite covering that is required, in terms
of the associated Weyl group and its monodromy action.
Much later this was extended to good characteristics [SB01]. 
All of this depends upon knowing that the formal, or
{\'e}tale, equivalence class of the singularity is determined
by the combinatorial structure of the exceptional locus;
that is, the singularities are \emph{taut}.
However, in bad characteristics they are not taut
(although Artin, following work of Lipman [Li69], showed
[Ar77], by giving a complete and explicit list of equations,
that there are only finitely many isomorphism
classes of each combinatorial type).

For types $A$ and $D$ Tyurina also showed
that simultaneous resolution could be achieved without
introducing algebraic groups into the picture by
manipulating explicit invariant polynomials under the Weyl group
and then making a suitable blow-up; this
relied on knowing the equation defining the singularity.
For type $E$ this approach has not been carried out; instead,
Brieskorn and Tyurina independently used the idea 
of embedding the singularity into a del Pezzo
surface. All of this too depends upon tautness. 
Each approach reveals that the necessary finite cover can be taken
to have the corresponding Weyl group as Galois group.

Artin's approach requires neither that the singularities
be taut, nor that they be realized on the unipotent variety
of a simply connected simple algebraic group, nor that they
be embeddable in a del Pezzo surface. In fact, there are du Val singularities
of type $E_8$ in characteristic $2$ which have none of these properties.
This is because there are five types of $E_8$ singularity, while only three
exist on a del Pezzo surface and localizing the unipotent variety of the
simple algebraic group $E_8$ at the generic point of the subregular locus
gives, of course, a unique singularity. (In fact,
although this plays no r{\^{o}}le here,
in any characteristic the unipotent
singularity of $E_8$ lies on a del Pezzo surface $D$ 
of degree $1$ where the $j$-invariant
of the anti-canonical curves is non-constant,
and this property, subject to the presence of an $E_8$ singularity on $D$, 
specifies $D$ uniquely. In turn, this pins down the 
singularity of the unipotent variety: in the notation of
[Ar77] it is $E_8^r$
where $r=4,2,1$ when the characteristic is $2,3,5$. For the details, see [GS].)

Rather, Artin's approach only depends upon the singularity being rational.
However, according to Artin, ``This more precise result
[concerning the Weyl group] does not follow directly
from our method''. The point of this paper is that
in fact Artin's methods do give this result.
The idea is merely to embellish the definition 
of Artin's functor $\Res$ so as to
include a suitable marking and polarization.

\begin{acknowledgments}
I am very grateful to Dave Benson, Ian Grojnowski,
Ian Leary and Michael Rapoport
for valuable correspondence and suggestions. 
\end{acknowledgments}
\end{section}
\begin{section}{Simultaneous resolution}\label{simul}
Fix a field $k$ and a surface $X_s/k$
with rational singularities.
Let $X'_s\to X_s$ denote the minimal resolution,
$F$ the exceptional locus in $X'_s$ and $E\subseteq F$
the maximal sublocus of $F$ on which $K_{X'_s}$ is trivial.
We shall assume that $E$ consists of copies of $\P^1_k$,
so that $E$ consists of the $(-2)$-curves $E_1,...,E_l$
in $X'_s$.
Of course, this assumption is vacuous if $k$ is algebraically closed.

Define the \emph{simple roots} to be the classes
of the $(-2)$-curves $E_1,...,E_l$ 
and set $Q=\oplus \Z E_i$,
the \emph{root lattice}.
Put $P=\Hom(Q,\Z)$, the \emph{weight lattice}.
Then 
there is a natural embedding
$\iota:P\inj\Pic_{X'_s}$ 
and the intersection pairing
identifies $Q$ with a sublattice of $P$.
There is a dual basis $\{\varpi_i\}$ of $P$ with 
$\deg\iota(\varpi_i)\vert_{E_j}=(\varpi_i,E_j)=\delta_{ij}$;
the $\varpi_i$ are the \emph{fundamental dominant weights}.

If $\Sigma$ is the set of simple roots then 
the reflexions in the elements of $\Sigma$
give a Coxeter system $(W,\Sigma)$
that acts on $P_\R$ and tesselates $P_\R$ into
chambers which are permuted simply transitively 
by $W$.

If $\sD$ is any such chamber,
$\G$ is the subgroup of the orthogonal group
$O_{P_\R}$ that preserves $P$ and $Q$
and $\Stab_{\sD}$ is the subgroup of $\G$
that preserves $\sD$, then
$\G$ is a semi-direct product $\G=W\rtimes\Stab_{\sD}$.

\def\tres{\widetilde{\Res}}
\def\RR{\tres_{\sD}}
\def\RRR{\widetilde{\Res}}

Now suppose that $S$ is 
a scheme of finite type
over $\L$,
that $s$ is a closed point of $S$ with $k(s)=k$, that
$f:X\to S$ is flat,
that $X_s$ is the fibre over $s$ 
and that the relative singular locus is
finite over $S$.

Fix a chamber $\sD$ as above and a
prime $\ell$ that is invertible in $\sO_S$.

There are various stacks over $S$ that we shall consider:
Artin's functor $\Res_{X/S}$;
$\Res_{P}$; $\Res_{P,\sD}$.
Here are their definitions.
\begin{enumerate}
\item
A $T$-point of $\Res_{X/S}$ is an 
isomorphism class of minimal resolutions
$$
\xymatrix{
{\widetilde{X_{T}}}  \ar[r]^{\pi}\ar[dr]_\Phi & 
{X_{T}}\ar[d]^\a\ar[r] & X\ar[d]\\
& {T}\ar[r] & {S}
};
$$
in particular, $\Phi$ is smooth, $\pi$ is projective
and is birational, in the sense that
$\pi_*\sO=\sO$, and $\pi$ induces a minimal resolution
of each geometric fibre of $\a$.

Given a $T$-point of $\Res_{X/S}$,
define $\NS_\ell({\widetilde{X_{T}}}/T)$
to be
the image of $\Pic({\widetilde{X_{T}}})$ in
$R^2\Phi_*\Z_\ell(1)$. This is a finitely generated $\Z$-module.

\item
A $T$-point of $\Res_{P}$ consists of a $T$-point of $\Res_{X/S}$
and a homomorphism
$\tphi:P\to \Pic_{\widetilde{X_{T}}}/\Pic_T$
such that the composite $\phi:P\to \NS({\widetilde{X_{T}}}/T)$
satisfies the following three conditions:

\noindent (i) $\phi(Q)$ lies in the image of
$R^2\Phi_c\Z_\ell(1)$;

\noindent (ii) $\phi(Q)$ is orthogonal to the relative canonical class
$K_{\widetilde{X_T}/T}$;

\noindent (iii) the composed pairing
$$P\times Q\to R^2\Phi_*\Z_\ell(1)\times R^2\Phi_c\Z_\ell(1)\to\Z_\ell$$
factors through $P\times Q\to\Z\to\Z_\ell$;

It is routine to write down the definition
of a morphism in each of these stacks, and to verify that
appropriate morphisms can be composed.

\item
Let $\Res_{P,\sD}$ denote the stack obtained from $\Res_P$ by adding
a fourth condition:

\noindent (iv) the cone $\langle \sD,K\rangle$ spanned by $\phi_\R(\sD)$
and the canonical class $K=K_{\widetilde{X_{T}}/T}$
lies in the nef cone $\NE(\widetilde{X_{T}}/X_T)$ of $\pi$.
\end{enumerate}

A $T$-point of $\Res_{P,\sD}$ is a family of $P$-marked,
$\langle\sD, K\rangle$-polarized surfaces and a $T$-point of $\Res_P$ is a 
family of $P$-marked surfaces.

For each $\sD$
there is a forgetful morphism $j_\sD:\Res_{P,\sD}\to\Res_P$, which is an open
immersion. There are also morphisms $q:\Res_P\to \Res_{X/S}$ and 
$r_\sD=q\circ j_\sD:\Res_{P,\sD}\to\Res_{X/S}$.

According to [Ar74] the stack
$\Res_{X/S}$ is represented by a locally quasi-separated
algebraic space $R$ over $S$ such that, for every field $K$,
$R\times\Sp K\to S\times\Sp K$ is an isomorphism.

Let $s$ also denote the unique point of $R$ over $s$.


Suppose that
$S^0\subset S$ is the complement of the discriminant
locus $\delta_S$ and $R^0$ its inverse image in $R$: then Artin
proves also that $R^0\to S^0$ is an isomorphism.

\begin{lemma}\label{nbhd}
Assume that $H^2(X'_s,\sO_{X'_s})=0$.
Then the forgetful morphism
$q:\Res_P\to\Res_{X/S}$
is a torsor under $\G$
over some Zariski neighbourhood
$U$ of $s$ in $S$.
\begin{proof}
The action of $\G$ on $\Res_P$
is given by the left action of $\G$
on the set
of homomorphisms $\phi:P\to\NS$
that satisfy conditions (i) and (ii):
$\g(\phi)=\phi\circ \g^{-1}$.
This makes it clear that $q$ is a pseudo-torsor
under $\G$.

The fibre of $q$ over $s\in S$ is non-empty
so it is enough to prove that
$q$ is dominant. For this
it is enough
to show that, given a henselian local
$S$-scheme $(T,s)$, every $T$-point of $\Res_{X/S}$
lies in the image of $q$.

Suppose given a minimal resolution
$\pi:\widetilde{X_T}\to X_T$.
The resolution $\tX_{s}\to X_s$
is minimal, so there is a unique homomorphism
$\tphi_0:P\to \Pic_{\tX_{s}}$ 
such that $\tphi(E_i)=E_i$.
(Recall that $Q$ is regarded as
a subgroup of $P$.)
Since $T$ is henselian,
the obstruction to extending
$\tphi_0(\varpi_i)$ to a class in 
$\NS(\widetilde{X_T}/T)$ lies
(after passing to formal completions and 
then algebraizing)
in $H^2(X'_{s},\sO_{X'_s})$, which vanishes. 
Therefore $\tphi_0$ extends to 
$\tphi:P\to \Pic_{\widetilde{X_T}}$.
\end{proof}
\end{lemma}
\begin{corollary}\label{indep} The definition of $\Res_P$ is independent of $\ell$.
\begin{proof} Suppose that $\ell'$ is another prime. 
Let $\Res_{P,\ell}$ and $\Res_{P,\ell'}$ be the 
corresponding functors.
Then define a third functor $\Res_P'$ as follows: an object
of $\Res_P'$ consists of a pair $({\widetilde{X_T}},\tphi)$
as in the definition of $\Res_P$ but demand that the induced
homomorphism
$$\phi':P\to\NS_\ell(\widetilde{X_T}/T)\oplus\NS_{\ell'}(\widetilde{X_T}/T)$$
satisfies
\begin{enumerate}
\item $\phi'(Q)$ lies in the image of $R^2\Phi_c\Z_\ell(1)\oplus R^2\Phi_c\Z_{\ell'}(1)$,
\item $\phi(Q)$ is orthogonal to $K_{{\widetilde{X_T}}/T}$ and
\item the composed pairing 
\begin{eqnarray*}
{P\times Q}&\to&{(R^2\Phi_*\Z_\ell(1)\oplus R^2\Phi_*\Z_{\ell'}(1))\times
(R^2\Phi_*\Z_\ell(1)\oplus R^2\Phi_*\Z_{\ell'}(1))}\\
&\to& {\Z_\ell\oplus\Z_{\ell'}}
\end{eqnarray*}
factors through $\Z\to \Z_\ell\oplus\Z_{\ell'}$.
\end{enumerate}
There are obvious forgetful morphisms
$\a:\Res_P'\to\Res_{P,\ell}$ and $\beta:\Res_P'\to\Res_{P,\ell'}$;
by Lemma \ref{nbhd} each of these morphisms
is a morphism of $\G$-torsors over $\Res_{X/S}$
and so is an isomorphism.
\end{proof}
\end{corollary}
{\bf{For the rest of this section
assume that the groups $H^2(X'_s,\sO_{X'_s})$ and
$H^2(X'_s,T_{X'_s}(-\log F))$
both vanish.}}
Of course, these assumptions hold if $X_s$ is affine,
or a partial resolution of an affine surface.

There is a unique maximal neighbourhood $U$
given by Lemma \ref{nbhd}.
From now on we replace $S$ by $U$,
so that $q$ is a $\G$-torsor over $S$.

In consequence of Lemma \ref{nbhd}, the forgetful morphisms
$\Res_P\stackrel{q}{\to}\Res_{X/S}$ and
$\Res_{P,\sD}\stackrel{r_\sD}{\to} \Res_{X/S}$ are {\'e}tale.
Since $\Res_{X/S}$ is represented by $R$,
$\Res_P$ and $\Res_{P,\sD}$ are then represented by
locally quasi-separated algebraic spaces
$R_P$ and $R_{P,\sD}$ over $S$, which are
{\'e}tale over $R$.
They form a diagram
$$\xymatrix{
{R_{P,\sD}}\ar@{^(->}[r]^{j_\sD}|{\circ}\ar[dr]_{r_\sD} & {R_P}\ar[d]^q\\
{} & {R}\ar[r] & {S.}
}$$

\begin{lemma}\label{separated}
\part[i] $R_{P,\sD}\to S$ is separated and quasi-finite.
\part[ii] $R_{P,\sD}$ is a scheme.
\part[iii] $R_P=\cup_\sD R_{P,\sD}$.
\begin{proof} In \DHrefpart{i} the separatedness is a consequence of the fact 
that a cone $\langle\sD,K\rangle$
of polarizations is specified in the data
that define a point of $\Res_{P,\sD}$.
The quasi-finiteness is immediate,
as is \DHrefpart{ii}.

\DHrefpart{iii} follows from the fact that
any resolution $\widetilde{X_T}\to X_T$
is a projective morphism.
\end{proof}
\end{lemma}

For any $T\to S$, let $T(s)$ denote the closed fibre over $s$
and let $S^h$ denote the henselization of $S$ at $s$.

\begin{lemma} \label{no1} 
Suppose that $Y$ is a quasi-finite
separated $S^h$-scheme of finite type.
Then

\part[i] there is
a unique minimal open subscheme $Y^1$ of $Y$ that contains $Y(s)$;
\part[ii] $Y^1$ is henselian, semi-local and finite over $S^h$;
\part[iii] $Y^1$
is also the maximal open subscheme of $Y$
that is finite over $S^h$;
\part[iv] $Y^1$ is the union of those connected
components of $Y$ that meet $Y(s)$.
\begin{proof} This is a consequence of Grothendieck's version of
Zariski's Main Theorem.
\end{proof}
\end{lemma}

Writing $(R_{P,\sD}\times_S S^h)^1$ as the disjoint union
of its connected components gives a decomposition
$$(R_{P,\sD}\times_S S^h)^1=\sqcup_{g\in\Stab_\sD}R^0_{P,\sD,g},$$
where each $R^0_{P,\sD,g}$ is a local henselian scheme that
is finite over $S^h$ and
$\Stab_\sD$ acts freely on $(R_{P,\sD}\times_S S^h)^1$
by permuting the various connected components $R^0_{P,\sD,g}$.

\begin{lemma}
$R_P\times_S S^h=\cup_\sD (R_{P,\sD}\times_S S^h)^1.$
\begin{proof}
$R_P\times_S S^h=\cup_\sD (R_{P,\sD}\times_S S^h)$,
so that $\cup_\sD (R_{P,\sD}\times_S S^h)^1$
is the unique minimal open subscheme of the non-separated scheme
$R_P\times_S S^h$ that contains $R_P(s)$,
which is a $\G$-torsor over a point.
Since $R_P\times_S S^h\to R\times_S S^h$ is a $\G$-torsor
the result follows.
\end{proof}
\end{lemma}

\begin{lemma} 
\part[i] If $T\to S^h$ is finite
then every point $T\to R$ factors through
$R^h$.

\part[ii] $R^h$ represents 
the restriction of $\Res_{X/S}$ to the category
of finite $S^h$-schemes.
\begin{proof} 
According to [Ar74], p. 332,
the co-ordinate ring $\sO_{R^h}$ is the direct limit
$\lim_{\to}\G(V,\sO_V)$ as $V$ runs over all affine
{\'e}tale neighbourhoods of $s$ in $R$.
So $R^h=(R\times_S S^h)^h$.
Since $R\to S$ is quasi-finite and surjective,
$R^h\to S^h$ is finite; 
therefore this system of rings
is cofinal with the subsystem consisting of
those rings
$\G(V,\sO_V)$ that are finite over $S^h$
and \DHrefpart{i} follows. 
\DHrefpart{ii} is an immediate consequence.
\end{proof}
\end{lemma}

In each chamber $\sD$, choose a fundamental domain
$\sS$ for the action of $\Stab_\sD$ on $\sD$
such that $\sD$ is tesselated by copies
of $\sS$ that are permuted simply transitively by $\Stab_\sD$
and the corresponding tesselation of $P_\R$
is preserved by $\G$. 
For example, take the case of $D_4$ and 
label the central vertex as $2$.
Then $\sD$ is the $\R$-span
of $\varpi_1,...,\varpi_4$
and we can take $\sS$ to be the cone spanned by
$\varpi_2,\varpi_1,\varpi_1+\varpi_4,\varpi_1+\varpi_3+\varpi_4$.

Then an \emph{orientation} of the nef cone $\sC$
of $X'_s\to X_s$ is the choice of fundamental
domain $\sS$ for the action of $\Stab_\sD$ on $\sD$
where $\sD$ is the positive chamber defined by
the irreducible $(-2)$-curves $E_1,...,E_r$ on $X'_s$
when these curves are regarded as simple roots.

Consider the functor $\Res_+$ defined on the category
of finite $S^h$-schemes $T$ as follows:
its $T$-points are elements of $R^h(T)$
together with an orientation of the nef cone of 
${\widetilde{X_T}}\to X_T$.

The natural morphism $\epsilon:\Res_+\to R^h$ is a torsor
under $\Stab_\sD$, so $\Res_+$ is represented by
a semi-local scheme $R^h_+=\sqcup_{g\in\Stab_\sD}R^h_g$,
a disjoint union of copies of $R^h$.

\begin{lemma} Each forgetful morphism
$p_\sD:(R_{P,\sD}\times_S S^h)^1\to R^h_+$ 
and each forgetful morphism
$p_{\sD,g}:R_{P,\sD,g}^0\to R^h_g$
is an isomorphism.
\begin{proof} From the definition of 
$R^h_+$ and the construction of $(R_{P,\sD}\times_S S^h)^1$
as an open subspace of $(R_{P,\sD}\times_S S^h)$
and from consideration of the functors that
they represent, it is clear that $p_\sD$
is {\'e}tale. Since both $(R_{P,\sD}\times_S S^h)^1$
and $R^h_+$  are semi-local
and $p_\sD$ is an isomorphism on the fibers over $s$
the result follows for $p_\sD$.
Since $p_{\sD,g}$ 
is the restriction of $p_\sD$ to connected components
the lemma is proved for it too.
\end{proof}
\end{lemma}

By construction, there is a family of
$P$-marked, $\sD$-polarized
surfaces over $R_{P,\sD}\times_S S^h$
whose fibers are resolutions of the fibers of $X\times_S S^h\to S^h$;
via the previous lemma there is then, for each
chamber $\sD$, a family $X'_\sD\to R^h$
of $P$-marked, $\sD$-polarized surfaces
whose closed fiber is $X'_s$.

Note that, for any chambers $\sD,\sE$, 
the isomorphisms $p_\sD$
and $p_\sE$  agree
on the overlap $(R_{P,\sD}\times_S S^h)^1\cap (R_{P,\sE}\times_S S^h)^1$
and so glue to a morphism 
$$p:\cup_\sD (R_{P,\sD}\times_S S^h)^1=R_P\times_S S^h\to R^h_+.$$
Given an algebraic space $Y$,
we say that a morphism $Y\to Q$ of algebraic spaces is the
\emph{separated quotient} of $Y$ if $Q$ is separated
and every morphism from $Y$ to a separated
algebraic space factors uniquely through $Q$.

\begin{lemma}\label{quotient}
$p:R_P\times_S S^h\to R^h_+$ is {\'e}tale and is
the separated quotient of $R_P\times_S S^h$.
\begin{proof}
This follows from the fact that for each chart $(R_{P,\sD}\times_S S^h)^1$
of $R_P\times_S S^h$ the morphism $(R_{P,\sD}\times_S S^h)^1\to R^h_+$
is an isomorphism.
\end{proof}
\end{lemma}

\begin{remark} Even when it exists, 
the quotient morphism to
a separated quotient is not always quasi-finite.
For example, the separated quotient of
the universal family over $R_P\times_S S^h$
is the pullback of the given family over $S$.
\end{remark}

We gather this information into a diagram
$$\xymatrix{
{(R_{P,\sD}\times_S S^h)^1}\ar@{^(->}[r]|-{\circ}\ar[drr]^{p_\sD}_{\cong}&
{\cup_{\sD} (R_{P,\sD}\times_S S^h)^1}\ar[r]^-{=}& 
{R_P\times_S S^h}\ar[d]^p\ar[dr]\ar[r]&{R_P}\ar[dr]^q\\
{}&{} & {R^h_+}\ar[r]^{\epsilon}&{R^h}\ar[r]& {R}\ar[r]&{S}
}$$
where $R^h$ is a henselian
local scheme, 
$(R_{P,\sD}\times_S S^h)^1$ and $R^h_+$ are henselian semi-local schemes, 
each finite and {\'e}tale over $R^h$,
$R_P\times_S S^h$ is a non-separated
scheme, $R_P$ and $R$
are locally quasi-separated
algebraic spaces, $q$ is a torsor under $\G$,
$\epsilon$ is a torsor under $\Stab_\sD$
and $p$ is a separated quotient.

\begin{lemma}\label{quotient2}
\part[i] $R^h_+$ is the normalization of
$R^h$ in $R_P\times_S S^h$. If $S^h$ is normal then
$R^h_+$ is also the normalization of $S^h$ in $R_P\times_S S^h$.
\part[ii] $\G$ acts on $R^h_+$.
\begin{proof} 
\DHrefpart{i} is clear from Lemma \ref{quotient}.
\DHrefpart{ii} follows from \DHrefpart{i}
and the fact that $\G$ acts on $R_P$.
\end{proof}
\end{lemma}

Now suppose that $X\to S$ is versal at $s$,
with respect to deformations over $\Lambda$.
Assume also that
\begin{enumerate}
\item either $X_s$ is affine
\item or the miniversal deformation space 
$Def_{X_s}$ of $X_s$ is formally smooth over $\Lambda$.
Note that this latter condition holds if $X_s$ is affine with only RDPs.
\end{enumerate}
Then [Ar74] the image of $R\times_S S^h$ in $S^h$
is an irreducible component $A$ of $S^h$,
and $A=S^h$ in the second case.
Let $\tA\to A$ denote the normalization. 
Then all spaces appearing in the diagram above, 
except maybe $S$, are smooth
in the appropriate sense.

\begin{theorem}\label{Main}
\part[i]\label{main} There is an effective action of $W$ 
on $R^h$ such that the geometric quotient
$[R^h/W]$ is naturally isomorphic to $\tA$.
\part[ii]\label{useful} $R^h$ is smooth over $\Lambda$ and is
the base of a versal deformation of the minimal resolution $X_s'$. 
\part[iii]\label{subsidiary} 
For every positive root $r$
with corresponding reflexion $\s_r$,
the fixed locus of $\s_r$ is the 
locus $D_r^h$ in $R^h$ that consists of
deformations of $X'_s$
where the root $r$ survives as an effective curve.
This locus $D_r^h$ is a divisor
and is smooth over $\Lambda$.
If $\{r_1,...,r_l\}$ is the set of simple roots
corresponding to the chamber $\sD$ then
the divisors $D_{r_1},...,D_{r_l}$
are transverse relative to $\Lambda$.
\part[iv]\label{part3} If $Def_{X_s}$ is formally smooth 
over $\Lambda$ then $S^h=[R^h/W]$.
\begin{proof}
For \DHrefpart{i}, observe first that 
the $W$-action on $R^h$ arises from
the $\G$-action on $R^h_+$ and the fact that
$R^h_+\to R^h$ is a torsor under $\Stab_\sD$,
so that the normal subgroup $W$ of $\G$
acts on each connected component of $R^h_+$.
The effectivity of the action follows from the fact
that $\G$ acts freely on the complement of the discriminant
in $R^h_+$.

Since $R\to A$ is an isomorphism over the complement of the discriminant,
it follows that $\deg((R_{P,\sD}\times_S S^h)^1\to A)=\#\G$.
So $\deg(R^h\to A)=\#W$, and then
$[R^h/W]\to \tA$ is an isomorphism by Galois theory.

For \DHrefpart{ii}
we start by copying
the proof of Lemma 3.3 of [Ar74].

Consider the 
formal deformation
functor for $X_s'$.
This has a hull $\Def_{X'_s}$, which
is formally smooth since, from the short exact sequence
$$0\to T_{X'_s}(-\log F)\to T_{X'_s}\to\oplus N_{F_j/X'_s}\to 0$$
and our assumption that $H^2(X'_s,T_{X'_s}(-\log F))$ vanishes,
the obstruction
space $H^2(X'_s,T_{X'_s})$ also vanishes.
Let ${\widehat{R_P}}$
denote the completion of $R_P\times_S S^h$ at any one of 
its closed points $x$.
Since $R_P\times_S S^h\to R^h$ is an isomorphism
in an {\'e}tale neighbourhood of $x$,
${\widehat{R_P}}$ is identified with
the completion of $R^h$.        
There are morphisms 
$$
\xymatrix{
{\Def_{X_s'}}& {\widehat{R_P}}\ar[l]_-{\hb}\ar[r] & {\widehat{S}}
}$$      
where $\widehat{S}$ is the completion of $S$ at $s$ and
$\hb$ is provided by
the semi-universal property of a hull. By
Lemma 3.3 of [Ar74] $\hb$ is formally smooth
and \DHrefpart{ii} is proved.

In $\Def_{X'_s}$ there is a divisor $\hD_i'$, 
formally smooth over $\Lambda$,
which is the locus where the exceptional $(-2)$-curve $E_i$ survives.
The Zariski tangent space to $\hD_i'$ is
$H^1(X'_s,T_{X_s'}(-\log E_i))$
and the Zariski tangent space to the locus $\cap \hD_i'$
where each $E_i$ survives is
$H^1(X'_s,T_{X_s'}(-\log\sum E_i)).$

Put $\hD_i=\hb^{-1}(\hD_i')$.
\begin{lemma}\part[i]
The natural homomorphism
$$H^1(X'_s,T_{X_s'}(-\log\sum E_i))
\to H^1(X'_s,T_{X_s'})$$
is injective and its image is of codimension $l$.
\part[ii] $\hD_1,...,\hD_l$ are transverse
divisors in ${\widehat{R_P}}$.
\begin{proof}
Take the cohomology of the exact sequence
$$0\to T_{X'_s}(-\log E)\to T_{X'_s}\to\oplus N_{E_i/X'_s}\to 0.$$
The lemma follows from the facts that
$\sN_{E_i/X'_s}\cong\sO_{\P^1}(-2)$
and, since, by assumption, the space $H^2(X'_s,T_{X'_s}(-\log F)$
is zero, $H^2(X'_s,T_{X'_s}(-\log E)$ is also zero.
\end{proof}
\end{lemma}

The irreducible curves are 
simple roots in $Q$.
They define a chamber $\sD$.
We shall identify $R_{P,\sD,1}^0$ with $R^h$.

Let $\s_i$ denote the reflexion in the simple root $\a_i$.

In $R_P$ there is an effective divisor $D_{P,i}$
defined by
$$D_{P,i}(T)=\{(\widetilde{X_T},\tphi)\vert\tphi(\a_i)
\ \textrm{is the class of an effective divisor}\}.$$
Then $\hD_i=D_{P,i}\times_{R_P}\widehat{R_P}$.
So $D_{P,i}$ is smooth over the coefficient ring $\Lambda$
and the divisors $D_{P,1},...,D_{P,l}$ are transverse.
Let $D_{+,i}$ denote the image of $D_{P,i}\times_S S^h$ in $R^h_+$
and $D_i^h$ its image in $R^h$.
By construction, $D_i^h$ is the locus in $R^h$
that parametrizes the deformations of $X'_s$
where the simple root $r_i$ survives as an effective
curve.
 
We have observed that $p:R_P\times_S S^h\to R^h_+$ is 
{\'e}tale and that $R^h_+$ is the largest
separated quotient of $R_P\times_S S^h$,
so that $p$ is $\G$-equivariant.
Moreover, $p$ is an isomorphism over the complement
of the discriminant in $S^h$ (or $A$).

Therefore, given $r\in R^h_+$ and $\g\in\G-\{1\}$,
$r\in\Fix_\g$ if and only if there is a 
henselian trait $V=\{0,v\}$
and two morphisms $a,b:V\to R_P\times_S S^h$ 
such that $a\ne b$, $a(v)=b(v)$
and $a(0)=b(0)\circ\g$.

We show next that the generic point $r_{+,i}$
of $D_{+,i}$ lies in 
the fixed locus $Fix_{+,\s_i}$ of $\s_i$ acting
on $R^h_+$.
Since $\G$ acts freely on $R_P\times_S S^h$, the point
$r_{+,i}$ lies in $Fix_{+,\s_i}$ if,
by the remark just made,
there are
two morphisms $a,b:V\to R_P\times_S S^h$
such that $a\ne b$, $p\circ a=p\circ b$,
$r=p\circ a(0)=p\circ b(0)$
and $p\circ a(v)\ne p\circ b(v)$.

To find these morphisms, take any such $V$
and a morphism $f:V\to R^h_+$ such that
$f(0)=r_{+,i}$ and
otherwise $f$ is in general position;
then the closed fibre of $X_V\to V$
has a single $A_1$ singularity and the
generic fibre is smooth.
It is well known 
(``the existence of flops'')
that if a smoothing of an $A_1$
singularity possesses a resolution, then it has
two such, and they are not isomorphic.
Therefore $f:V\to R^h_+$ has two liftings
$a,b:V\to R_P\times_S S^h$ as described above,
so that $r_{+,i}\in Fix_{+,\s_i}$.
It follows that $D_i^h$ is contained in
the fixed point locus $Fix_{\s_i}^h$ of
$\s_i$ acting on $R^h$.

To complete the proof of \DHrefpart{iii}, we need to show that
$\Fix_{\s_i}^h\subset D_i^h$.

Suppose that $\eta=\eta_1=[E_1]$
is a simple root. Suppose also that
$x\in (R_{P,\sD}\times_S S^h)^1$ with $x\mapsto t\in S$,
where $\sD=\sum_{i=1}^l\R_{\ge 0}\varpi_i$.
Assume that $\eta$ is not the class of a $(-2)$-curve
on the minimal resolution $X'_t$. That is,
$E_1$ does not survive as an effective class when
$X'_s$ deforms to $X'_t$.

Now $\sD$ is contained in the nef cone 
$\NE(X'_t)$, by the definition of $(R_{P,\sD}\times_S S^h)^1$. 

\begin{lemma}
$\s_{\eta}(\sD)\subseteq\NE(X'_t)$.
\begin{proof}
$\s_{\eta}(\sD)$ is spanned by
$\varpi_1+\eta_1,\varpi_2,...,\varpi_l$.
Any irreducible $(-2)$-curve $F_t$ on $X'_t$ specializes
an effective cycle $F_s$ on $X'_s$; then $[F_s]$ is
a positive root $\phi$ and $\phi\ne\eta_1$.
So $\phi.\varpi_j\ge 0$ for all $j$, 
while $\phi.(\varpi_1+\eta_1)\ge 0-1=-1$
since $\phi\ne\eta_1$. Suppose 
$\phi.(\varpi_1+\eta_1)<0$; then
$\phi.\varpi_1=0$ and $\phi.\eta_1=-1$.
But $\phi=\sum_{j\ge 1} n_j\eta_j$ with $n_j\ge 0$,
so that $\phi=\sum_{j\ge 2} n_j\eta_j$ and 
$\phi.\eta_1=-1$, which is absurd.
\end{proof}
\end{lemma}

Therefore $\s_{\eta}(x)\in (R_{P,\sD}\times_S S^h)^1$ also;
since $\G$ acts freely on $\cup_{\sD}(R_{P,\sD}\times_S S^h)^1=R_P\times_S S^h$
it follows that $\s_{\eta}(x)\ne x$.
However, the morphism $p_{\sD}:(R_{P,\sD}\times_S S^h)^1\to R^h_+$
is an isomorphism, so $\s_{\eta}(p(x))\ne p(x)$.
Hence $Fix_{\s_i}^h\subseteq D_i^h$,
so that $Fix_{\s_i}^h= D_i^h$.

So part \DHrefpart{iii} of the theorem
is proved for every simple reflexion.

Every reflexion $\s_r$ in a positive
root $r$ is conjugate in $W$ to
a simple reflexion, 
so the locus $Fix_{\s_r}^h$
is a smooth divisor $\tD^h_r$.
In characteristic zero we know, by Proposition 2.4 (i)
of [Wa79], that $\tD^h_r=D^h_r$;
since $D^h_r$ is also a smooth divisor, it follows that
$D^h_r=Fix_{\s_r}^h$, and \DHrefpart{iii} is proved.

Finally, if $Def_{X_s}$
is formally smooth over $\Lambda$ then $S^h=\tA$,
and now Theorem \ref{Main} is proved.
\end{proof}
\end{theorem}

\begin{corollary}\label{2.12} Suppose given a normal local
henselian scheme $(T,0)$ and a family $g:Y\to T$
of surfaces such that the closed fibre
$Y_0$ has only du Val singularities
and the generic fibre is smooth.
Then $Y\to T$ has a resolution if and only if
the Galois action on $H^2$ of the generic
fibre is trivial.
\begin{proof} The question is local
on $Y$, so we may assume that $Y\to T$ is pulled back
via a morphism $T\to S$ whose image does not lie
in the discriminant locus of $f:X\to S$. On the
one hand, the
$W$-covering $R_\sD^0\to S^0$ is exactly
the covering defined by the monodromy action
on $H^2$ of the geometric generic fibre of
$X\to S$, while on the other $Y\to T$ has a resolution
if and only if the map $T\to S$ factors through
$R_\sD$. Since $T$ is normal, we are done.
\end{proof}
\end{corollary}

\begin{remark}
\noindent (1) Note that for families that map to the discriminant
locus in $S$, it might be
necessary to take an inseparable cover; for example, this
happens for the family $xy+z^2+t=0$
of $A_1$ singularities in characteristic $2$.

\noindent (2) Suppose that $X\to C$ is a morphism from
the germ of a smooth threefold
to the germ of a smooth curve, that the closed fiber $X_s$ has a du Val
singularity and that $\ch k=0$. Then the monodromy (the image of a generator
of the local $\pi_1$, the fundamental group of the punctured curve)
is a Coxeter element of $W$ [Dm75].
However, in positive or
mixed characteristic the local $\pi_1$ is a local Galois group
and is not cyclic and it is not clear how to describe the image of this
Galois group in $W$. For example, if the residue characteristic is $3$ and
the type of the singularity is $A_2$ then the image of Galois equals $W$.
\end{remark}
\end{section}
\begin{section}{Some topology of the situation}
Suppose in this section that $X\to S$ is versal at $s$,
that $H^2(X'_s,T_{X'_s}(-\log E))$ and $H^2(X'_s,\sO_{X'_s})$
are both trivial, that $X_s$ has only RDPs and that $S$ is a
local henselian scheme.

We have seen that the non-separated scheme
$R_P$ is obtained by
glueing copies $R_{P,\sD}^1$ of the (separated) scheme $R^h_+$ 
along open subschemes. 
Moreover, 
\begin{enumerate}
\item $W$ acts on each connected component
$R_{P,\sD,g}^0$ of $R_{P,\sD}^1$,
where $g$ runs over $\Stab_{\sD}$,
and the connected components 
$R_{P,\sD,g}$ are all isomorphic, say to 
$R_{P,\sD}^0$,

\item the image of $R_{P,\sD,g}^0$ in $R^h_+$
is a connected component $R^h_g$ of $R^h_+$,

\item $R^h_g$ 
is isomorphic to $R^h$ and is the maximal
separated quotient of a connected component $R_P^0$
of $R_P$ and

\item $W$ acts on $R^h$.
\end{enumerate}
Now that we know, thanks to Theorem \ref{Main},
that $R^h$ contains a hyperplane arrangement
consisting of the divisors $D^h_r$
parametrized by the positive roots $r$, 
we can describe this glueing more precisely,
as an analogue of the construction of
the real algebraic prevariety $\sZ(\sA)$ of [Pr07].

\begin{proposition}\label{nonsep}
$R_P^0$ is the non-separated scheme
obtained by glueing copies $R^h_\sD$ of $R^h$,
one copy for each chamber $\sD$ in the Euclidean vector
space $P_\R$, as follows.

For chambers $\sD,\sE$ in $P_\R$,
the intersection $R^h_{\sD}\cap R^h_{\sE}$
is given by
$$R^h_{\sD}\cap R^h_{\sE}=R^h-\cup_r D^h_{r}$$
where $r$ runs over those positive roots $r$ such that in
$P_\R$ the wall $H_r$
separates $\sD$ and $\sE$.
\begin{proof} It's enough to prove the analogous statement
for the possibly disconnected schemes 
$(R_{P,\sD})^1$.

Observe that $(R_{P,\sD})^1\cap (R_{P,\sE})^1$
consists of those points $x$ in $R_P$
that map to a point $t\in S$ such that on the minimal
resolution $X'_t$ of $X_t$ the chambers
$\sD$ and $\sE$ both lie in $\NE(X'_t)$.
This is equivalent to saying that
$\sD$ and $\sE$ both lie on the positive
side of the wall $H_r$ if the positive root
$r$ corresponds to an effective cycle on $X'_t$,
and there is nothing left to prove.
\end{proof}
\end{proposition}

\begin{corollary} $W$ acts freely on $R^0_P$
and $R=R^0_P/W$.
\noproof
\end{corollary}

Now suppose that $k=\R$ and that $X_s$
has a du Val singularity. According
to Brieskorn, Slodowy \emph{et al.}
we can, by abuse of notation, write $S=[\mathfrak t/W]$ 
and $R^h=\mathfrak t$, where $\mathfrak t$
is a Cartan subalgebra of the
Lie algebra of the relevant split simple algebraic group. 
So $R^h$ is a complexified
hyperplane arrangement. Let $S_0\subseteq S$
be the complement of the discriminant
and $\tS_0$ the inverse image of $S_0$ in $R^h$.

\begin{corollary}\label{2.11}
$S_0(\C)$ is weakly homotopy
equivalent to $R(\R)$.
\begin{proof} By Corollary \ref{nonsep},
$R_P^0(\R)$ is nothing but the non-separated
manifold $\sZ(\sA)(\R)$.
The main result of [Pr07] is that
there is a $W$-equivariant map 
$\tS_0(\C)\to\sZ(\sA)(\R)$ that is a weak homotopy equivalence.
Since $W$ acts freely on both sides,
taking quotients by $W$ gives a weak homotopy equivalence
$$S_0(\C)\to [\sZ(\sA)(\R)/W]=[R_P^0(\R)/W]= R(\R).$$
\end{proof}
\end{corollary}
\begin{corollary} $S_0(\C)$ is a $K(\pi,1)$
where $\pi$ is the corresponding generalized braid group.
\begin{proof}
This follows from Deligne's result [D72]
that $R(\R)$ is a $K(\pi,1)$.
\end{proof}
\end{corollary}
\end{section}
\begin{section}{Polynomial rings of $W$-invariants over $\Z$ via RDPs}\label{dp}
Consider the polynomials $f$ in $\Z[x,y,z]$ given 
in the attached table.
\begin{table*}[h!] 
\begin{center} 
\begin{tabular}{cccc}
{Type}&{$f$}&
{$\Pi_{fund}$}&{$\Pi_{extra}$}\\
{$A_n$}&{$xy+z^{n+1}$}&
{$(2,...,n+1)$}&{$(1)$}\\
{$D_{2n}$}&{$x^2+z^2y+zy^n$}&
{$(2,4,...,4n-2,2n)$}&{$(1^{(2)},3,5,...,2n-1)$}\\
{$D_{2n+1}$}&{$x^2+z^2y+y^nx$}&
{$(2,4,6,...,4n-2,4n,2n+1)$}&{$(1,3,5,...,2n-1)$}\\
{$E_6$}&{$x^2+xz^2+y^3$}&
{$(2,5,6,8,9,12)$}&{$(1,2,3,4,6)$}\\
{$E_7$}&{$x^2+z^3y+y^3$}&
{$(2,6,8,10,12,14,18)$}&{$(1,3,4,5,9)$}\\
{$E_8$}& {$x^2+y^3+z^5$}&
{$(2,8,12,14,18,20,24,30)$}&{$(3,4,5,6,9,10,15)$}
\end{tabular}
\end{center}
\end{table*}
In each case the surface $X_0$ defined by $f=0$ in $\A^3_\Z=\Sp\Z[x,y,z]$
has an effective action of $\GG_{m,\Z}=\Sp\Z[\lambda^\pm]$ 
and has an RDP of the indicated type at every field-valued
point of the origin $0$ in $\A^3_\Z$.
Further calculation shows that, in each case, there is, for some $N$, a
$\GG_{m,\Z}$-equivariant family $\sX\to\A^N_\Z$ that is versal
(taking $\Z$ to be the coefficient ring) at every
field-valued point of the origin $0_S\cong\Sp\Z$ in $\A^N_\Z$
and that
\begin{enumerate}
\item no weight of the action on $\A^N_\Z$ is zero and
\item the set $\Pi$ of positive weights of this action is 
$\Pi=\Pi_{fund}\cup\Pi_{extra}$, counting multiplicities,
where $\Pi_{fund}$ and $\Pi_{extra}$ are as tabulated.
\end{enumerate}
$\Pi=\Pi_{fund}$
is the set of fundamental degrees
(exponents plus $1$) of the corresponding root system
[GrLie4-6] and $\Pi_{extra}$ is the set of extra weights.

Write $N=N^++N^-$ where $N^+=\#\Pi$ is the number of positive weights
and $\rho=N^+-r=\#\Pi_{extra}$.

Because having RDPs and versality are both open conditions,
there is a $\GG_{m,\Z}$-invariant
open subscheme $S$ of $\A^{N}_\Z$ such that $S$ contains $0_S$,
the induced family $X\to S$ is everywhere versal and all its geometric
members are affine surfaces with RDPs.

Since none of the weights listed above is zero,
the origin $0_S$ is the fixed locus of the $\GG_{m,\Z}$-action
on $\A^{N}_\Z$, and so of the action on $S$.

Fix a prime $\ell$. Then, over $S[1/\ell]=S\otimes\Z[1/\ell]$
there is a non-separated scheme $R_P=R_P^{(\ell)}$
given by the construction of Section \ref{simul}.

Note that $\GG_{m,\Z[1/\ell]}\times \G$ acts
on $R_P^{(\ell)}$ and the morphism $R_P^{(\ell)}\to S[1/\ell]$
is $\GG_{m,\Z[1/\ell]}$-equivariant.

Let $\pi^{(\ell)}:Q^{(\ell)}\to S[1/\ell]$ 
denote the normalization of $S[1/\ell]$ in $R_P^{(\ell)}$.
The next lemma is well known.
\define\bX{\overline{X}}
\begin{lemma}\label{normal}
If $G$ is a smooth affine group scheme over some affine 
normal base $\Sp B$
and $X\to Y$ is a dominant quasi-finite
$G$-equivariant morphism of affine normal $G$-schemes
over $B$
then $G$ acts on the normalization $N$ of $Y$ in $X$.
\begin{proof}
$\sO_{N}=\{f\in\sO_X\vert f\ \textrm{is\ integral\ over}\ \sO_Y\}$.
Let $f\in\sO_{N}$
and suppose that $\sum_0^r a_if^i=0$
with $a_i\in\sO_Y$ and $a_r=1$. Let $\mu_Z:\sO_Z\to\sO_Z\otimes_B\sO_G$
be the co-actions, for $Z=X,Y$. Then
$\sum\mu_Y(a_i)\mu_X(f)^i=0$,
so that $\mu_X(f)$ lies in $\sO_{X}\otimes\sO_G$
and is integral over $\sO_Y\otimes\sO_G$.
That is, $\mu_X(f)$ is in the normalization
of $\sO_Y\otimes\sO_G$ in $\sO_X\otimes\sO_G$.
But this normalization is $\sO_{N}\otimes\sO_G$,
since $G$ is smooth over $B$, and therefore
$\mu_X$ restricts to a co-action
on $\sO_{N}$.
\end{proof}
\end{lemma}

\begin{proposition}\label{glueable}
\part[i] $Q^{(\ell)}$ is smooth over $\Z[1/\ell]$.
\part[ii] $Q^{(\ell)}$ is the separated quotient
of $R_P^{(\ell)}$.
\part[iii] $\GG_{m,\Z[1/\ell]}\times \G$ acts on $Q^{(\ell)}$ and
$\pi^{(\ell)}:Q^{(\ell)}\to S[1/\ell]$ is $\GG_{m,\Z[1/\ell]}$-equivariant.
\part[iv] $\G$ acts effectively on $Q^{(\ell)}$ and
$\pi^{(\ell)}:Q^{(\ell)}\to S[1/\ell]$ identifies $S[1/\ell]=[Q^{(\ell)}/\G]$.
\begin{proof} 
According to Lemmas \ref{quotient}
and \ref{quotient2}, $R^{(\ell),h}_+$ is,
for all henselizations $S^h$ of $S[1/\ell]$,
the normalization
of both $S^h$ and $R^{(\ell),h}$ in $R_P^{(\ell)}\times_{S[1/\ell]} S^h$.
Now $R^{(\ell),h}_+\to R^{(\ell),h}$ is a $\Stab_\sD$-torsor, so {\'e}tale,
and $R^{(\ell),h}$ is smooth, and so $Q^{(\ell)}\times_{S[1/\ell]} S^h$ is smooth.
\DHrefpart{i} follows.

\DHrefpart{ii} can be checked after passing to $S^h$,
where it follows from Lemma \ref{quotient}

For \DHrefpart{iii}, it is clear that $\G$ acts. 
The existence of the $\GG_{m,\Z[1/\ell]}$-action
follows from the definition of $Q^{(\ell)}$ 
as the normalization of $S[1/\ell]$ in $R^{(\ell)}_P$.

\DHrefpart{iv} follows from the facts that $\pi^{(\ell)}$ is finite,
its degree is $\#\G$ and $S[1/\ell]$ is smooth, so normal.
\end{proof}
\end{proposition}

Now suppose that $\ell'$ is a second prime.
Then over $S[1/\ell\ell']$ the schemes
$R_P^{(\ell)}$ and $R_P^{(\ell')}$ are isomorphic,
by Corollary \ref{indep}.
Therefore $Q^{(\ell)}$ and $Q^{(\ell')}$ are isomorphic
over $S[1/\ell\ell']$ and so can be glued to give
$\pi:Q\to S$.

\begin{proposition}\label{glued}
\part[i] $Q$ is smooth over $\Z$.
\part[ii] $Q$ is the separated quotient
of $R_P$.
\part[iii] $\GG_{m,\Z}\times \G$ acts on $Q$ and
$\pi:Q\to S$ is $\GG_{m,\Z}$-equivariant.
\part[iv] $\G$ acts effectively on $Q$
and $\pi:Q\to S$ identifies $S=[Q/\G]$.
\begin{proof} This is an immediate consequence of
Proposition \ref{glueable}.
\end{proof}
\end{proposition}

Since the weights of the $\GG_{m,\Z}$-action on $\A^{N}_{\Z}$
are never zero
the fixed locus of the action is $0_S$.
\define\nbd{neighbourhood}

Since $S$ is a $\GG_{m,\Z}$-equivariant \nbd\ of 
$0_S$ in $\A^{N}_{\Z}$, it follows that
$S$ contains $S^+$ and $S^-$,
where 
$$S^\pm=\{x\in\A^{N}_{\Z}\vert
\lim_{\lambda^\pm\to 0}\lambda(x)\subset 0\}.$$
Then $S^+$ (resp., $S^-$) is defined 
as a subscheme of $S$ by the vanishing of the co-ordinates
of negative (resp., positive) weight, so 
$S^\pm\cong\A^{N^\pm}_{\Z}$
and $S^+$ and $S^-$ intersect transversely in the origin $0_S$.

\begin{lemma}\label{saturated} $\pi^{-1}(0_S)_{red}\cong\Sp\Z\times\Stab_\sD$.
\begin{proof} Calculation shows that,
in each case, the minimal
resolution $X'_0\to X_0$ is obtained
by successively blowing-up along
copies of $\Sp\Z$, so that $X'_0$
is smooth over $\Z$.
Then
for any chamber $\sD$, the set of
markings $\phi:L\to\NS(X'_0)$ 
such that $\phi(\langle \g,\sD\rangle)$ lies 
in the nef cone (and then equals the nef cone)
is a torsor under $\Stab_\sD$.
So $X'_0$ defines a $\Sp\Z\times\Stab_\sD$-point of $R_{P,\sD}$.
Similarly,
for any field $k$ and for every morphism
$f:\Sp k\to 0_S$, $X_0\otimes k$ has a unique
minimal resolution. Therefore the set of lifts of $f$
to $\pi^{-1}(0_S)$ is a torsor under $\Stab_\sD$, 
and therefore $\pi^{-1}(0_S)_{red}\to\Sp\Z\times\Stab_\sD$
is an isomorphism. 
\end{proof} 
\end{lemma}

Define $Q^\pm=\pi^{-1}(S^\pm)$
and $0_Q=\pi^{-1}(0_S)_{red}$.
Choose a connected component $\tQ$ of $Q$,
and write $\tQ^\pm=\tQ\cap Q^\pm$
and $0_{\tQ}=\tQ\cap 0_Q$.
Then $\GG_{m,\Z}\times W$ acts on $\tQ$ and on $\tQ^\pm$.
Observe that 
$$\tQ^\pm =\{x\in \tQ\vert\lim_{\lambda^\pm\to 0}\lambda(x)\subset 0_{\tQ}\},$$
so that
the closure of each $\GG_{m,\Z}$-orbit
in $\tQ^\pm$ meets $0_{\tQ}$.

\begin{lemma} The fixed locus of
$\GG_{m,\Z}$ acting on $\tQ$
is $0_{\tQ}$.
\begin{proof} $0_S$ is the fixed locus of
$\GG_{m,\Z}$ acting on $S$ and the fixed locus of
$\GG_{m,\Z}$ on $Q$ is smooth over $\Sp\Z$.
\end{proof}
\end{lemma}

The next lemma is a version of Theorem 2.5 of [KR82],
but in mixed characteristic.

\begin{lemma}\label{KR} 
Suppose that $\GG_{m,\Z}$ acts on a smooth affine
$\Z$-scheme $X=\Sp A$ such that the fixed locus
of the $\GG_{m,\Z}$-action is isomorphic to $\Sp\Z$
and meets the closure of every orbit.
Then $X$ is $\GG_{m,\Z}$-equivariantly isomorphic
to an affine space over $\Z$.
\begin{proof}
The hypotheses imply that $A$
is a graded ring, say $A=\oplus_{n\ge 0}A_n$,
that $A_0=\Z$ and that the ideal $I$ of 
$Fix(\GG_{m,\Z}\vert_{X})$
is $I=\oplus_{n\ge 1}A_n$.
Now the argument from [KR82] goes through
to show that $X$ is isomorphic to an affine
space over $\Sp\Z$,
since their Corollary 1.4 is stated and proved
over any base.
\end{proof}
\end{lemma}

\begin{lemma} $\tQ^+\cap \tQ^-=0_{\tQ}$ as schemes
and $\tQ^\pm$ is isomorphic to $\A^{N^\pm}_{\Z}$. 
\begin{proof}
Certainly the Krull dimension
of $\tQ^\pm$ is $\dim \tQ^\pm=\dim S^+=N^\pm +1$.

From the description of $\tQ^\pm$ as limit loci
it follows that the tangent $\Z$-module
$T_{0_{\tQ}}\tQ^\pm$ is the part of the representation $T_{0_{\tQ}}\tQ$
of $\GG_{m,\Z}$ where the weights are all of the corresponding sign.
No weight of this tangent action is zero,
and so $\tQ^+,\tQ^-$ intersect transversely in $0_{\tQ}$
and are therefore
smooth along $0_{\tQ}$ of the dimensions already indicated.

Since $\tQ^\pm$ is smooth along $0_{\tQ}$ and on each one of them
$0_{\tQ}$ meets the closure of any given $\GG_{m,\Z}$-orbit
it follows from Lemma \ref{KR} that they are affine spaces
over $\Z$.
\end{proof}
\end{lemma}

Since $W$ acts on the henselization $Q_t^h$
of $Q$ at each field-valued point $t$ of $\Phi$
as a Coxeter system $(W,\Sigma=\{\s_1,...,\s_8\})$
(that is, the fixed locus of each $\sigma\in\Sigma$
acting on $Q_t^h$
is a smooth divisor and these divisors are transverse),
it follows that $(W,\Sigma)$ also acts as a Coxeter
system on $\tQ^+$. 

\begin{lemma}
For any field $k$ the action of $W$ on $\tQ^+\otimes k$
is effective.
\begin{proof} 
The construction so far has been made for deformations
where the coefficient ring is $\Z$.
When repeated with $\Z$ replaced by $k$
as coefficient ring,
that is, when we consider only deformations
over $k$ of $X_0\otimes k$,
the spaces $R_P, Q$ etc., are replaced
by $R_P\otimes k, Q\otimes k$ etc.
Then the result follows from 
the effectivity of the action of $\Stab_\sD$ given by
Proposition \ref{glued}.
\end{proof}
\end{lemma}

Let $M$ denote a root lattice of type $ADE$
and $r$ its rank. 
In the polynomial ring $\Z[M]$
regard the elements of $M$ as being of degree $1$
and so write $\Z[M]=\Z[1^r]$.
\begin{theorem}\label{extend}
The action of $W$ on the root lattice $M$
extends to a graded action of $W$ on a polynomial ring
$\sO=\Z[M][\Pi_{extra}]=\Z[1^r\cup\Pi_{extra}]$ on $N^+$ variables
over $\Z$ such that

\part[i] $\sO^W$ is a polynomial $\Z$-algebra 
$\Z[\Pi_{fund}\cup\Pi_{extra}]$ and
\part[ii] \label{extend ii}
for all normal domains $A$, $W$ acts effectively on $\sO\otimes A$
and $(\sO\otimes A)^W=\sO^W\otimes A$.
\begin{proof}
We deduce this from three lemmas.
\begin{lemma}\label{4.6} $S^+=[\tQ^+/W]$
and $S^+\otimes_\Z A=[(\tQ^+\otimes_\Z A)/W]$
for all normal domains $A$.
\begin{proof} Since $W$ acts effectively on $\tQ^+$
the commutative diagram
$$\xymatrix{
{\tQ^+}\ar[r]\ar[d]&{Q}\ar[d]^-{\pi}\\
{[\tQ^+/W]}\ar[r]&{[Q/\G]=S}
}$$
is Cartesian in a \nbd\ of the generic point of $[\tQ^+/W]$
and so the morphism $[\tQ^+/W]\to S$ identifies $[\tQ^+/W]$
with the normalization of its image, which is $S^+$.

The same argument applies after tensoring with $A$.
\end{proof}
\end{lemma}

Write $\sO_{\tQ^+}=\Z[x_1,...,x_{N^+}]$.
This is a positively graded polynomial $\Z$-algebra 
with a graded $W$-action, where $\deg x_i\ge 1$ for all $i$.
For each $\sigma_i\in\Sigma$ and for every henselization $Q^h_t$
at a point $t$ of $\Phi$,
the fixed locus $Fix(\s_i\vert_{Q^h_t})$ contains $(Q^-)^h_t$,
as already remarked. So $Fix(\s_i\vert_Q)$ is a smooth divisor
that contains $Q^-$
and the divisors $Fix(\s_i\vert_Q)$, for $i=1,...,r$,
are transverse. 
Set $D_{\s_i}=Fix(\s_i\vert_{\tQ^+})$;
then $D_{\s_1},...,D_{\s_r}$
are smooth transverse divisors in $\tQ^+$.
Say $D_\s=(f_\s)_0$; then
each $f_\s$ is unique up to an element of
$\sO_{\tQ^+}^*=\pm 1$,
so that for each $\s$ there is a character
$\chi_\s$ of $\GG_{m,\Z}$ such that
$\lambda(f_\s)=\chi_\s(\lambda)f_\s$.
That is, each $f_\s$ is homogeneous.
Since the divisors $D_{\s_i}$ are transverse,
we can then assume $x_i=f_{\s_i}$ for $i=1,...,8$.
The reflexions in $W$ are conjugate, so
all of $x_1,...,x_r$ have equal degrees.
Put $M_1=\sum_1^r\Z x_i$.

\begin{lemma} $W$ acts 
on the subring $\Z[x_1,...,x_r]$ of $\sO_{\tQ^+}$
via its standard action on the root lattice $M$.
\begin{proof} The fixed locus $Fix(W\vert_{\tQ^+})$ is 
given by
$Fix(W\vert_{\tQ^+})=\cap_1^r D_{\s_i}=\cap_\s D_\s$,
where in the second intersection
$\s$ runs over all reflexions in $W$.
So for every reflexion $\s$,
$f_\s$ lies in the ideal $M_1.\sO_{\tQ^+}$;
since $\deg f_\s=\deg x_1$, it follows that
$f_\s=\sum\a_ix_i$ and every $\a_i$ lies in $\Z$.
Since $Fix(w\s w^{-1})=w(Fix(\s))$
it follows that $M_1$ 
is a representation of $W$.
Since the $x_i$ are algebraically independent
they are certainly linearly independent,
so that $\rk M_1=r$.

Moreover, $\s_i$ preserves and acts trivially on the hyperplane
$(x_i)_0$ in $\tQ^+$, so that
for every $g\in\sO_{Q^+}$
$\s_i(g)-g$ lies in the ideal $(x_i)$.
In particular,
$\s_i(x_j)=x_j+b_{ji}x_i$ for some $b_{ji}\in\Z$ 
if $i\ne j$ and $\s_i(x_i)=-x_i$.

If $i\ne j$ and $(\s_i\s_j)^2=1$ then
$b_{ij}=b_{ji}=0$, 
while if $(\s_i\s_j)^3=1$ then
$b_{ij}b_{ji}=1$, so that $b_{ij}=b_{ji}=\pm 1$.
Since the Coxeter diagram is a tree
we can choose $b_{ij}=1$ for all such $i,j$,
and then deduce that
the representations $M_1$ and $M$
of $W$ are isomorphic.
\end{proof}
\end{lemma}\label{degree}
\begin{lemma} 
$\deg x_i=1$ and
$\sO_{\tQ^+}=\Z[M][y_1,...,y_\rho]$
where the degrees of the extra generators
$y_1,...,y_\rho$ run over $\Pi_{extra}$.
\begin{proof}
The degrees can be read by tensoring with $\Q$.
\end{proof}
\end{lemma}

The theorem now follows from Lemmas \ref{4.6}
and \ref{degree}
by taking $\sO=\sO_{\tQ^+}$.
\end{proof}
\end{theorem}

%
%
%
%
\end{section}
\begin{section}{Chow rings of complete flag varieties}
It is well known that the structure of rings of
invariants under Weyl groups
is connected with the topology of the corresponding 
simply connected split reductive group.
For example, if $M$ is either the root or 
the weight lattice of a split semisimple
algebraic group $G$ over any base, then $\Q[M]^W$
is isomorphic to $A^*(G/B)\otimes\Q$,
where $A^*$ denotes the Chow ring.
(Recall that, for flag varieties, 
$H^*((G/B)(\C),\Z)$ is naturally isomorphic
to the Chow ring, up to a doubling of degrees.
Therefore we do not need to be concerned here with the base
over which $G$ is defined.)

Given an $\N$-graded ring $A=\oplus_{n\in\N}A_n$,
set $A_+=\oplus_{n>0}A_n$.
If a group $\G$ acts on $A$
set $A_\G=A/I$, where $I=A.A^\G_+$.
This is the ring of \emph{co-invariants}.

Demazure [Dm73] constructed
a finite graded $\Z$-algebra $H=H(\Z[M])$
and a graded $W$-homomorphism $c:\Z[M]\to H$,
all purely in terms of $M$ and $W$.
He then showed that $H$ is isomorphic to 
the integral cohomology ring $H^*(G/B,\Z)=A^*(G/B)$
and that
$c$ induces an injective homomorphism
$\bar{c}:\Z[M]_W\inj H$
whose image is the subring generated by
the first Chern classes of line bundles on $G/B$
that can carry $G$-linearizations.

As is well known, $c$ (or $\bar{c}$)
is not surjective for many groups;
this is related to the fact that the ring
$\Z[M]^W$ of invariants is not a polynomial ring,
even when $M$ is the weight lattice.
From now on take $M$ to be the root lattice.
In Section \ref{dp} we have constructed
a polynomial $\Z[M]$-algebra $\sO$ such that
$\sO^W$ is a polynomial $\Z$-algebra;
what we prove here is that $c$
extends to a surjective homomorphism
$c:\sO\to H$ and that $\bar{c}:\sO_W\to H$
is an isomorphism.

\begin{theorem}\label{abstract} Suppose that $M$
is the root lattice belonging to
a finite Coxeter system $(W,\Sigma)$ of type $ADE$
and that $\sO$ is a graded polynomial $\Z[M]$-algebra
on $\rho$ further generators and that
$W$ acts on $\sO$ compatibly with its action
on $\Z[M]$. Let $\pi:\Sp\sO\to\Sp\Z[M]$
be the projection and assume that
\begin{enumerate}
\item\label{hyp 1} for every $\s\in\Sigma$, the fixed locus
$Fix(\s\vert_{\Sp\sO})$ equals
the inverse image
$\pi^{-1}(Fix(\s\vert_{\Sp\Z[M]}))$,

\item
$W$ acts effectively
on $\Sp(\sO\otimes R)$ for every normal domain
$R$ and 

\item $\sO^W$ is a polynomial
$\Z$-algebra.
\end{enumerate}
Then there is a surjective $W$-equivariant homomorphism
$c:\sO\to H$ 
that induces an isomorphism $\bar{c}:\sO_W\to H$.
\begin{proof}
We begin by recapitulating results from [Dm73]
and extending some of them slightly to cover our situation.
These extensions are easy.

\begin{lemma} 
$\sO_W$ is flat over $\Z$.
\begin{proof} $I=\sO.\sO^W_+$ defines the subscheme
$\pi^{-1}(0)$ of $\Sp\sO$. By
Lemma \ref{saturated} this is $1$-dimensional and is quasi-finite over $\Sp\Z$;
since $I$ is generated by $r+\rho$ elements these generators form a
regular sequence, and the result follows.
\end{proof}
\end{lemma}

Say $N$ is the number of reflexions in $W$.
This is also the number of positive roots in $M$.
The fixed locus $Fix(\s\vert_{\Sp\sO})$
equals the zero locus $(x_\s)_0$
of some element $x_\s\in M\subset\sO_1$.
Define $d=\prod_\s x_\s\in\sO_N$ and $J=\sum\det(w)w$.

Since any reflexion $\s$ acts trivially on the divisor $(x_\s)_0$
in $\Sp\sO$,
$u-\s(u)$ lies in the ideal $(x_\s)$, so that there is
an $\sO^W$-linear divided difference operator $\Delta_\s:\sO\to\sO$ 
of degree $-1$ defined,
up to $\pm 1$, by
$$\Delta_\s(u)=(u-\s(u))/x_\s.$$
If $\s=\s_i\in\Sigma$ then write $\Delta_{\s_i}=\Delta_i$.





The formulae $(3)-(6)$ of [Dm73]
for the maps $\Delta_\s$ are valid.

let $\sD$ denote the $\sO^W$-algebra
of $\sO^W$-linear endomorphisms of $\sO$
generated by $\sO$ and the $\Delta_\s$.
It is a left $\sO$-module and is graded, since
each $\Delta_\s$ is homogeneous
of degree $-1$.

\begin{lemma}\label{Lemme 2} (= Lemme 2 of [Dm73]) For all $\Delta\in\sD$
there exist $(\Delta'_i,\Delta''_i)_{i=1,...,n}$ such that
$\Delta(uv)=\sum\Delta'_i(u)\Delta''_i(v).$
\noproof
\end{lemma}

Let $\epsilon:\sO\to\Z$ be the augmentation map.
Then $\epsilon\sD$ is a $\Z$-module of linear
maps $\sO\to\Z$ that kill $I$, so is finitely generated.
Say $H(\sO)$ is its $\Z$-dual. There is a dual
$\Z$-linear map $c=c_\sO:\sO\to H(\sO)$.
By Lemma \ref{Lemme 2} $H(\sO)$ is a graded co-algebra.

\begin{proposition} $H(\sO)$
has a unique structure as a graded commutative ring
for which $c$ is a graded ring homomorphism.
Moreover, $W$ acts on $H(\sO)$ in such a way
that $c$ is $W$-equivariant.
\begin{proof} Prop. 2 of [Dm73] and the remarks following.
\end{proof}
\end{proposition}

According to Th. 1 of [Dm73], there is a well defined
operator $D_w$ for each $w\in W$ such that $D_\s=\Delta_\s$
for each reflexion $\s$.
If $w_0$ is the longest element of $W$ then $D_{w_0}=J/d$
([Dm73], Prop. 3 (b)).

\begin{proposition} (= Cor. 4 of [Dm73])
\part[i] $\{\epsilon D_w\}_{w\in W}$ is a $\Z$-basis
of $\epsilon\sD$.
\part[ii] Write $z_w=\epsilon D_w$. Then  $\{z_w\}_{w\in W}$
is the unique $\Z$-basis 
of $H(\sO)$
such that $c$ is given by
$c(u)=\sum_{w\in W}\epsilon D_w(u)z_w.$
\part[iii] $H(\sO)_i$ is based by $\{z_w\vert\ell(w)=i\}$.
\noproof
\end{proposition}


\begin{proposition}\label{5.15}
\part[i] (Poincar{\'e} duality)
The multiplication $H(\sO)_i\times H(\sO)_{N-i}\to H(\sO)_N=\Z.z_{w_0}$ 
is a perfect pairing of $\Z$-modules.

\part[ii] The ring $H(\sO)$ is the same
for all polynomial $\Z[M]$-algebras $\sO$ that
satisfy Hypothesis \ref{hyp 1} of Theorem \ref{abstract}.
That is, there is a $W$-equivariant commutative diagram
$$\xymatrix{
{\Z[M]}\ar[r]^{c_{\Z[M]}}\ar[d]&{H(\Z[M])}\ar[d]^{\cong}\\
{\sO}\ar[r]_{c_\sO}&{H(\sO)}
}$$
where $H(\Z[M])\to H(\sO)$ is an isomorphism.
\begin{proof}
Cor., p. 293 of [Dm73], and Remarque 2) following.
\end{proof}
\end{proposition}
Let $H$ denote the common value of the rings $H(\sO)$ and
$H(\Z[M])$.
As explained in [Dm73], $H$ is isomorphic (after a doubling
of degrees) as a graded $\Z$-algebra with $W$-action to the
integral cohomology ring $H^*(G/B,\Z)$
where $G$ is the split simple algebraic group
of the given type $A,D$ or $E$.
So $\bar{c}$ is an injective graded $W$-equivariant homomorphism
$$\bar{c}:\sO_W\inj H^*(G/B,\Z)$$
of finite flat $\Z$-algebras.

The next three results are taken from or inspired by [Bro09].

\begin{lemma}\label{splitting} There is an $\sO^W$-linear splitting
$\nu:\sO\to\sO^W$ of the inclusion $\sO^W\inj\sO$
whose kernel is a free $\sO^W$-module.
\begin{proof}
Consider first the induced homomorphism
$j:\Z=\sO^W/\sO^W_+\to\sO_W$.
Via the augmentation map $\sO\to\Z$
with kernel $\sO_+$, $j$ has an $\sO^W$-linear
splitting. So there is a $\Z$-basis
$\{1,x_1,...,x_{\#W-1}\}$ of $\sO_W$
where each $x_i$ is homogeneous and lies in ${\sO_W}_+$.
Since $\sO^W\to\sO$ is finite,
any lifting of this $\Z$-basis 
to a subset $\Phi$ consisting of 
homogeneous elements of $\sO$
is, by the graded version of Nakayama's Lemma,
a generating set of $\sO$ as an $\sO^W$-module.
Now $\sO^W\to\sO$ is also flat,
since both rings are regular of the same dimension,
and therefore $\Phi$ is an $\sO^W$-basis of $\sO$.
In particular, there is an $\sO^W$-basis of $\sO$
that includes the element $1$. The lemma follows.
\end{proof}
\end{lemma}

\begin{lemma}\label{diff} 
The different $\mathfrak D$ of the morphism $\Sp\sO\to\Sp\sO^W$
is defined by the principal ideal $(d)$.
\begin{proof} By assumption, $W$ acts effectively
on $\Sp(\sO\otimes k)$ for all fields $k$, so that, in particular,
$\mathfrak D$ does not contain the mod $2$ fibre $\Sp(\sO\otimes\F_2)$.
Now on one hand $\mathfrak D$ is an effective $W$-invariant
Cartier divisor on $\Sp\sO$, since $\Sp\sO\to\Sp\sO^W$
is a finite flat and separable morphism of regular schemes,
and on the other hand the ideal $\sI_{\mathfrak D}$ of $\mathfrak D$ 
coincides with $(d)$
over $\Sp\Z[1/2]$, since the only elements of
$W$ that act on $\Sp\sO$ with fixed points in codimension one 
are the reflexions $\s$,
the fixed locus of each of which is the corresponding wall $(x_\s)_0$.
Therefore $\sI_\mathfrak D=(d)$ over $\Sp\Z$.
\end{proof}
\end{lemma}

\begin{proposition}\label{exists} There exists $a\in\sO_N$ such that $J(a)=d$.
\begin{proof} Write $S=\sO$, $R=\sO^W$.
Since $R\to S$ is finite and flat, the relative dualizing module
$\omega_{S/R}$ is the graded $S$-module 
$$\omega_{S/R}=\Hom_R(S,R).$$
Since $R$ and $S$ are smooth $\Z$-algebras,
Lemma \ref{diff} gives, by the well known canonical
isomorphism $\mathfrak D^{-1}\buildrel{\cong}\over{\to}\omega_{S/R}$,
an isomorphism
$$\phi:S.d^{-1}\buildrel{\cong}\over{\to}\omega_{S/R}$$
of graded $S$-modules defined by the $R$-bilinear pairing 
$$S.d^{-1}\times S\to R:(a/d,s)\mapsto \Tr(as/d)=\sum_{w\in W}w(as/d).$$
Choose $a\in S$ such that $\phi(a/d)=\nu$,
where $\nu$ is the splitting provided by Lemma \ref{splitting};
since $\nu(1)=1$, we get
$\Tr(a/d)=1$. Since $\phi$ is a graded isomorphism
we can choose $a$ to be homogeneous, and then 
$a\in S_N$. Since $w(d)=\det w.d$
the result follows from the definition of $J$.
\end{proof}
\end{proposition}

Since $\sO^W$ is a polynomial ring and $\sO$ is flat over $\sO^W$, 
the ideal $I$
is generated by a regular sequence in $\sO$,
so that $\sO_W$ is a finite flat complete intersection
$\Z$-algebra. The ring $H$ is, because multiplication
gives a perfect pairing into $H_N$, by Proposition \ref{5.15},
a finite flat Gorenstein $\Z$-algebra.

Now we can prove the theorem.
It is enough to show that
the injective homomorphism $\bar{c}:\sO_W\to H$ is 
surjective. 

By Proposition \ref{exists},
the fact that $D_{w_0}=J/d$
and the fact that,
by the definition of $H$, $H_N$ is generated
as a $\Z$-module by $D_{w_0}$, the map 
$\bar{c}_N:(\sO_W)_N\to H_N$
is surjective, 
and so is an isomorphism, since both sides
are torsion-free $\Z$-modules of rank one.
Therefore $\bar{c}:\sO_W\to H$ is a graded homomorphism
of finite flat graded Gorenstein $\Z$-algebras
such that, for every field $k$, the induced
homomorphism $\bar{c}\otimes 1_k:(\sO_W)\otimes k\to H\otimes k$
of finite local Gorenstein $k$-algebras
induces an isomorphism of socles.
Therefore $\bar{c}\otimes 1_k$ is injective.
Since both algebras have the same dimension over $k$
(namely, $\# W$), $\bar{c}\otimes 1_k$ is an isomorphism,
and we are done.
\end{proof}
\end{theorem}

\begin{corollary}\label{5.11}
There is an action of $W$ on a polynomial $\Z[M]$-algebra
$\sO=\Z[M][\Pi_{extra}]$ 
such that $\sO^W$
is a polynomial $\Z$-algebra
$\Z[\Pi_{fund}\cup\Pi_{extra}]$ 
and $H^*(G/B,\Z)$ is isomorphic to $\sO_W$. 
%
\begin{proof}
According to Theorem \ref{extend}, the hypotheses of Theorem \ref{abstract}
are satisfied, and the Corollary follows.
\end{proof} 
\end{corollary}
\end{section}
\begin{section}{Partial resolutions of RDPs and partial flag varieties}
A \emph{partial resolution} of an RDP over a field $k$ is a quasi-projective surface
$Y$ with RDPs that is a partial resolution of an affine surface $X$
that is defined over $k$ and has RDPs.
That is, there are proper birational morphisms
$$X'\to Y\to X$$
where $X$ is affine with RDPs and $X'$ is the minimal resolution
of both $X$ and $Y$.

If $X\to S$ is a flat family of affine surfaces where each
geometric fibre has RDPs then a partial resolution of $X$
over $S$
consists of $Y\to X$ where $Y\to S$ is flat and for each field-valued
point $s$ of $S$ the morphism $Y_s\to X_s$ is
a partial resolution.

Fix a coefficient ring $\Lambda$ for $k$.

\begin{lemma} The local deformation space of
a partial resolution $Y$ of an RDP over a field $k$ 
is formally smooth over $\Lambda$.
\begin{proof} $H^2(Y,\sF)=0$ for all coherent
sheaves $\sF$ on $Y$.
\end{proof}
\end{lemma}

Now take $\Lambda=\Z$ and
suppose that $X_0$ is one of the RDPs
over $\Z$
considered in Section \ref{dp}.
Then there are $2^r$ partial resolutions $Y$ of $X_0$
over $\Z$;
they correspond, via taking the curves
contracted by $X'\to Y$, to subsets $\Theta$ of the set of
simple roots. So, for example, $Y=X'$ if $\Theta=\emptyset$.

Suppose also
that $S$ is the base of the deformation
considered there,
that $S^h$ is the henselization of $S$ at a point $s$ of
$0_S$ and that
$T^h$ is the henselian base of a versal deformation of $Y_s$.
So, by the previous lemma, 
$S^h$ and $T^h$ are both smooth, in the henselian sense.

Suppose that
$W_\Theta\subseteq W$ is the subgroup of $W$ generated by the reflexions
in the members of $\Theta$.
By Theorem \ref{main}, $S^h\cong [R^h/W]$ 
and $T^h\cong [R^h_1/W_\Theta]$
where each of $R^h$ and $R^h_1$ is the base
of a versal deformation of $X'$.
Both $R^h$ and $R_1^h$ are smooth
since $H^2(X',\sF)=0$ for all coherent sheaves $\sF$ on $X'$.

For the rest of this section we use the notation of Section \ref{dp}.
\begin{theorem}\label{partial}
\part[i] $[R^h/W_\Theta]$ is smooth.
\part[ii]\label{smooth} $[\tQ/W_\Theta]$ is smooth over $\Sp\Z$.
\part[iii] $[\tQ^+/W_\Theta]$ is the spectrum of a polynomial
ring over $\Z$.
\begin{proof} The preceding discussion, combined with Theorem \ref{main},
shows that $[R_1^h/W_\Theta]$ is smooth. Then
\DHrefpart{i} follows from the fact that, since any two versal deformation spaces
of the same object are smoothly equivalent, 
$R^h_1$ is smoothly equivalent to $R^h$.
\DHrefpart{ii} follows from \DHrefpart{i} and the fact that
at every point lying over the origin $0_S$ of $S$ the henselization of $Q$
is smoothly equivalent to $R^h$. \DHrefpart{iii} is then
proved using a $\GG_{m,\Z}$-action, as in Section \ref{dp}.
\end{proof}
\end{theorem}
Fix a Borel subgroup $B$ of $G$ and 
let $P$ denote the parabolic subgroup of $G$ that contains $B$
and corresponds to $\Theta$
(so that, for example, $P=B$ if $\Theta=\emptyset$.)
Put $\sO=\sO_{\tQ^+}$.

\begin{corollary}\label{BGG} 
\part[i] $\sO^{W_\Theta}$ is a polynomial $\Z$-algebra.
\part[ii] $\sO^{W_\Theta}/\sO^W_+.\sO^{W_\Theta}$
is isomorphic to $H^*(G/P,\Z)$.
\begin{proof}
\DHrefpart{i} is just Theorem \ref{smooth} above.

For \DHrefpart{ii},
write $C=\sO,B=\sO^{W_\Theta}$ and $A=\sO^W$.
Then, since $H^*(G/P,\Z)=H^*(G/B,\Z)^{W_\Theta}$ ([BGG73], Theorem 5.5),
 there is a commutative diagram of rings
$$\xymatrix{
{C}\ar@{->>}[r]&{C/A_+.C}\ar[r]^<<<<<{\bar c}_<<<<<{\cong}&{H^*(G/B,\Z)}\\
{B}\ar@{^(->}[u]\ar@{->>}[r]&{B/A_+.B}\ar@{^(->}[u]\ar[r]^<<<<<{\bar b}&{H^*(G/P,\Z)}\ar@{^(->}[u]\\
{A}\ar@{^(->}[u]\ar@{->>}[r]&{A/A_+}\ar@{^(->}[u]\ar[r]^{\bar a}_{\cong}&{\Z.}\ar@{^(->}[u]
}$$

\begin{lemma} If $B_1$ is a $B$-algebra that is torsion-free
as a $\Z$-module, then $B_1=(C\otimes_B B_1)^{W_\Theta}$.
\begin{proof} By Lemma \ref{splitting}
there is a $B$-linear splitting $\nu:C\to B$
whose kernel $N$ is a free $B$-module.
So $C=B\oplus N$, and then $C\otimes_B B_1=B_1\oplus(N\otimes_B B_1)$.
Therefore the inclusion $B_1\inj C\otimes_B B_1$ is cotorsion-free
as a $\Z$-module, and then $B_1\inj (C\otimes_B B_1)^{W_\Theta}$
is also cotorsion-free
as a $\Z$-module. Now 
$$B_1\otimes_\Z\Q=(C\otimes_B B_1\otimes_\Z\Q)^{W_\Theta}=
(C\otimes_B B_1)^{W_\Theta}\otimes\Q,$$
since finite groups are linearly reductive over $\Q$,
and so $B_1=(C\otimes_B B_1)^{W_\Theta}$.
\end{proof}
\end{lemma}
Take $B_1=B/A_+.B$. Then
\begin{eqnarray*}
B/A_+.B&=&(C\otimes_B(B/A_+.B))^{W_\Theta}=(C/A_+.C)^{W_\Theta}\\
&=&
H^*(G/B,\Z)^{W_\Theta}=H^*(G/P,\Z).
\end{eqnarray*}
\end{proof}
\end{corollary}
\end{section}
\begin{section}{Moduli of Enriques surfaces}
Suppose that $k$ is an algebraically closed field
of characteristic $2$ and that $Y$ is a smooth Enriques
surface over $k$. Then $\NS(Y)$ is isomorphic to the even
unimodular lattice $E=E_{10}(-1)$. Denote by $O(E)$
the orthogonal group of $E$ and $O^+(E)$ the index $2$ subgroup of
$O(E)$ consisting of elements that preserve the two
cones of positive vectors in $E_\R$. It is known [CD89]
that $O^+(E)$ is the Weyl group $W(E)$, the group
generated by reflexions in the roots of $E$. 
Fix, once and for all, a chamber $\sD_0$ defined by the roots
(that is, the $(-2)$-vectors)
in the positive cone of $E\otimes\R$. (We shall not always be scrupulous in
distinguishing between $\sD$ and its closure.)  This defines a root basis
$\beta_1,\ldots,\beta_{10}$ of $E$ that in turn defines a Dynkin diagram of
type $E_{10}=T_{2,3,7}$. We recover $\sD_0$ as $\sD_0=\sum\R_{\ge 0}\varpi_i$, where
$\varpi_1,\ldots,\varpi_{10}$ are the fundamental dominant weights defined by
the root basis. That is, $\varpi_i.\beta_j=\delta_{ij}$.
We label the simple roots $\beta_1,\ldots,\beta_{10}$ according
to the diagram
{\input Dynkin
\setdynkin tags=y \setdynkin slope="25
\def\x#1{$\scriptstyle#1$}
\def\present#1{\hbox{\hbox to2cm{$#1$:\hfil}\hskip1.5cm}}
\par\medskip
\dynkinE 10 \x{\beta_1} \x{\beta_3} \x{\beta_4} \x{\beta_2} \x{\beta_5} \x{\beta_6}
\x{\beta_7} \x{\beta_8} \x{\beta_9} \x{\beta_{10}}
\par
\bigskip
}
According to [L15] and [EHS],
the weight $\varpi_1$ defines a Cossec--Verra polarization $\lambda$ on $Y$;
in turn this defines a birational contraction $Y\to Z_0$,
where $Z_0$ has du Val singularities, such that $\lambda$ descends
to an ample class $\lambda$ on $Z_0$.
The forgetful map $\Def_{Z_0,\lambda}\to\Def_{Z_0}$ 
of infinitesimal deformation functors is an
isomorphic, as is the corresponding morphism of
henselian functors,
and [Li10] these functors 
are formally smooth over either $\Sp\W(k)[[f,g]]/(fg-2)$ or $\Sp\W(k)$
according to whether $\Pic^\tau_{Y}$ is isomorphic to $\a_2$ or not.

Suppose that $Z\to S$ is a versal deformation of $(Z_0,\lambda)$
and that $S$ is henselian; then
consider the stack $\Res_\sD$ over $S$ that is defined much as before: the objects
over an $S$-scheme $T$ consist of a resolution
$\pi:{\widetilde{Z_T}}\to Z_T$ together with a choice
of chamber $\sD$ in the nef cone of ${\widetilde{Z_T}}\to T$
such that $\pi^*\lambda$ is the vector corresponding to $\varpi_1$
under the unique isomorphism $\sD_0\to \sD$.
As before, $\Res_\sD$ is represented 
by a finite local $S$-scheme $R_\sD$ that carries an action
of a Weyl group $W$ such that the geometric quotient $[R_\sD/W]$ is $S$.
The Weyl group is that associated to the configuration of
singularities on $Z_0$. The root lattice
corresponding to this configuration embeds into the root
lattice of type $D_9$, since $D_9$ is the orthogonal complement
of $\varpi_1$ in $E_{10}$.

The infinitesimal functors $\Def_{Y}$ and $\Def_{Y,\sD}$ 
are also isomorphic, so there is, as before, a diagram
$$
\xymatrix{
{\Def_{Y}} & {\widehat{R_{\sD}}}\ar[l]^{\a} \ar[d]\\
& {\widehat{S}}
}$$
where $\alpha$ is formally smooth.

We can summarize this discussion in the following result.

\begin{theorem}\label{enriques} 
A versal deformation space $\Def_Y$ of a smooth
Enriques surface $Y$ is smoothly equivalent to a complete
local scheme $R_\sD$ that carries an action of a Weyl
group $W$ such that $[R_\sD/W]$ is regular. Moreover,
$W$ belongs to a root sublattice of the root lattice
$D_9$ and $[R_\sD/W]$ is formally smooth over either $\W$ or
$\W[[f,g]]/(fg-2)$. 

The scheme $R_\sD$ is regular
unless $Y$ is both classical and exceptional
(i.e., $\Pic^\tau_Y\cong\Z/2$ and $Y$ has vector fields).
\begin{proof}
The only thing that has not been proved is that,
with the stated exceptions, $\Def_Y$ is regular.
This is proved in [EHS]. More precisely, if
$Y$ is an $\alpha_2$-surface, then a hull for
$\Def_Y$ is given by
$\W[[f,x_2,\ldots,x_{12}]]/(fg-2)$,
where $g$ is divisible by neither $2$ nor $f$,
and otherwise it is a power series ring over $\W$.
\end{proof}
\end{theorem}

\begin{remark} This picture describing the local moduli
of Enriques surfaces can now be described in terms
of the local picture above.

Suppose that $S_1$ is the base of a versal deformation
$Z_1\to S_1$ of a sufficiently small affine {\'e}tale neighbourhood of
$\Sing Z_0$.
Then there is a $W$-covering $R_{1,\sD}\to S_1$, where 
$R_{1,\sD}$ represents the stack of $\sD$-polarized resolutions
of $Z_1\to S_1$, and a commutative square
$$\xymatrix{
{R_\sD}\ar[r]\ar[d] & {R_{1,\sD}}\ar[d]\\
{S} \ar[r] & {S_1.}
}$$
Since both vertical arrows are $W$-coverings, the square is Cartesian
(which was not obvious \emph{a priori}). 
\end{remark}
\end{section}
\bigskip
\bibliography{alggeom,ekedahl}

\providecommand{\bysame}{\leavevmode\hbox to3em{\hrulefill}\thinspace}
\begin{thebibliography}{EGAIII:2}

\bibitem[Ar74]{Ar74}
M.~Artin, \emph{Algebraic construction of {B}rieskorn's resolutions}, J.
  Algebra (1974), \textbf{29}, 330--348.


\bibitem[Ar77]{Ar77}
\bysame, \emph{Coverings of the rational double points in characteristic $p$},
  Complex analysis \& algebraic geometry, Cambridge Univ. Press, Cambridge,
  1977.

\bibitem[BGG73]{BGG73}
I.N.~Bernstein, I.M.~Gel'fand and S.I.~Gel'fand, \emph{Schubert cells and 
cohomology of the spaces $G/P$}, Russian Math. Surveys \textbf{28} (1973), 1--26.

\bibitem[GrLie4-6]{BouL4-6}
N.~Bourbaki, \emph{Groupes et alg{\`e}bres de {L}ie. {IV}--{VI}}, Hermann,
  Paris, 1968.

\bibitem[Br68]{Br68}
E.~Brieskorn, \emph{Die Aufl{\"o}sung der rationalen Singularit{\"a}ten holomorpher Abbildungen}, Math. Ann. \textbf{178} (1968)
255--270.

\bibitem[Br70]{Br70}
\bysame, \emph{Singular elements of semi-simple algebraic groups}, Proc. ICM 1970, Nice.

\bibitem[Bro09]{Bro09}
A.~Broer, \emph{On Chevalley--Shephard--Todd's theorem in 
positive characteristic}, in Symmetery and Spaces, Birkh{\"a}user, 2009.

\bibitem[BR75]{BR75}
D.~Burns and M.~Rapoport, \emph{On the Torelli problem for k{\"a}hlerian K3 surfaces},
Ann. Sci. ENS \textbf{8} (1975), 235--273.


\bibitem[CD89]{CD89}
F.~R. Cossec and I.~V. Dolgachev, \emph{Enriques surfaces {I}}, Progress in
  {M}athematics, vol.~76, Birkhäuser, Boston Basel Berlin, 1989.

\bibitem[D72]{D72}
P.~Deligne, \emph{Immeubles des groupes de tresses g{\'e}n{\'e}ralis{\'e}s},
Inv. Math. \textbf{17} (1972), 273–-302.

\bibitem[Dm73]{Dm73}
M.~Demazure, \emph{Invariants sym{\'e}triques entiers des groupes de Weyl
et torsion}, Inv. Math. \textbf{21} (1973), 287--301.

\bibitem[Dm75]{Dm75}
\bysame, \emph{Classification des germes {\`a} point critique isol{\'e}
et {\`a} nombres de modules $0$ ou $1$ (d'apr{\`e}s V.I.~Arnol'd},
S{\'e}m. Bourbaki, LNM \textbf{431}, Springer, 1975.


\bibitem[DZ15]{DZ15}
H.~Duan and X.~Zhao, \emph{Schubert presentation of the cohomology ring of 
flag manifolds $G/T$},  LMS J. Comput. Math. \textbf{18} (2015), 489–506.




\bibitem[ESB99]{ESB99}
T.~Ekedahl and N. Shepherd-Barron, \emph{On exceptional {E}nriques
  surfaces}, arXiv:math/0405510 

\bibitem[EHS]{EHS}
T.~Ekedahl, M. Hyland and N. Shepherd-Barron, \emph
{Moduli and periods of simply connected Enriques surfaces},
arXiv:math':1210.0342





\bibitem[GS]{GS}
I.~Grojnowski and N.~Shepherd-Barron,
\emph{Del Pezzo surfaces as Springer fibres for exceptional groups}, arXiv:1507.01872 





\bibitem[L10]{L10}
C.~Liedtke, \emph{Moduli and lifting of Enriques surfaces},
 arXiv:1007.0787.

\bibitem[Li69]{Li69}
J.~Lipman, \emph{Rational singularities, with applications to algebraic
  surfaces and unique factorization}, Publ. Math. IHES \textbf{36} (1969),
  195--279.




\bibitem[Na10]{Na10}
M.~Nakagawa, \emph{The integral cohomology ring of $E_8/T$},
Proc. Japan Acad. \textbf{86} (2010), 64--68.





\bibitem[Pr07]{Pr07}
N.~Proudfoot, \emph{A non-Hausdorff model for the complement of a complexified 
hyperplane arrangement}, Proc. Amer. Math. Soc. \textbf{135} (2007), 3989--3994.





\bibitem[SB96]{SB96}
N.~I.~Shepherd-Barron, \emph{Some foliations on surfaces in characteristic
  $2$}, J. Algebraic Geom. \textbf{5} (1996), no.~3, 521--535.

\bibitem[SB01]{SB01}
\bysame, \emph{On simple groups and simple singularities},
Israel J. Math. \textbf{123} (2001), 179--188.

\bibitem[Sl80]{Sl80}
P.~Slodowy, \emph{Simple singularities and simple algebraic groups},
Lecture Notes in Math. \emph{815} (1980).


\bibitem[TW74]{TW74}
H.~Toda and T.~Watanabe, \emph{The integral cohomology ring of
$F_4/T$ and $E_6/T$}, J. Math. Kyoto Univ. \textbf{14} (1974), 257--286.


\bibitem[Ty70]{Ty70}
G.N.~Tyurina, \emph{ Resolution of singularities of flat
    deformations of double rational points}, 
Funkcional. Anal. i Prilozhen. \textbf{4} (1970), 77-–83.

\bibitem[Wa79]{Wa79}
J.M.~Wahl, \emph{Simultaneous resolution and discriminantal loci},
Duke Math. J. \textbf{46} (1979), 341--375.

\end{thebibliography}
\bibliographystyle{pretex}
\end{document}